\documentclass[11pt,reqno]{amsart}
\usepackage[shortlabels]{enumitem}
\usepackage{esint}
\usepackage{physics}
\usepackage{relsize}
\usepackage[backref]{hyperref}
\usepackage{mathtools}
\usepackage{amsfonts}
\usepackage[all]{xy}
\usepackage{geometry}
\usepackage{amsmath,amssymb} 
\usepackage{mathrsfs}
\usepackage{float}
\usepackage{faktor}
\usepackage{amscd}
\usepackage{tikz-cd}
\usepackage{multirow}

\newtheorem{theorem}{Theorem}[section]
\newtheorem{lemma}[theorem]{Lemma}

\theoremstyle{corollary}
\newtheorem{corollary}[theorem]{Corollary}
\theoremstyle{conjecture}
\newtheorem{conjecture}[theorem]{Conjecture}

\theoremstyle{assumption}

\theoremstyle{proposition}
\newtheorem{proposition}[theorem]{Proposition}
\theoremstyle{remark}
\newtheorem{remark}[theorem]{Remark}
\numberwithin{equation}{section}
\everymath{\displaystyle}

\newcommand{\ii}{\ensuremath{\sqrt{-1}}}
\newcommand{\pa}{\partial}
\newcommand{\pp}{\bar\partial}

\newcommand{\noop}[1]{}

\DeclareMathOperator{\Ric}{Ric}

\begin{document}

\title[Local Rigidity of Bergman Log-Flatness]{Local Rigidity of Bergman Logarithmic Flatness on Compact Hermitian Symmetric Spaces}

\author{Zhiqin Lu}
\address{Department of Mathematics, University of California, Irvine, Irvine, CA 92697, USA}
\email{zlu@uci.edu}

\author{Yuang Shi}
\address{School of Mathematical Sciences, Shanghai Jiao Tong University, Shanghai 200240, China}
\email{yuangshi@sjtu.edu.cn}

\begin{abstract}
We prove the local rigidity of Bergman logarithmic flatness near the standard metrics associated with fixed invariant positive polarizations on compact Hermitian symmetric spaces. Every sufficiently small smooth deformation preserving Bergman logarithmic flatness is induced by a holomorphic automorphism and a constant rescaling. In the Szeg\H{o} setting, local rigidity was established by Lu--Tian (2004) on $\mathbb{CP}^n$ and, more recently, by Loi--Palmieri--Zuddas (2026) on compact Hermitian symmetric spaces.
\end{abstract}

\maketitle
\tableofcontents

\section{Introduction}
\subsection{Bergman logarithmic flatness and rigidity}

The relation between Bergman kernels of positive line bundles and the K\"ahler geometry of the underlying manifolds has been
extensively studied; see, for example, \cite{Tian1990,Zelditch1998,Catlin1999,Lu2000,Donaldson2001,DaiLiuMa2006,BermanBerndtssonSjostrand2008,Xu2012,DonaldsonSun2014}. A related question concerns the extent to which the boundary behavior of the Bergman kernel of the associated disk bundle constrains the K\"ahler metric on the base. The logarithmic term in the Bergman kernel of this disk bundle provides a further link between the two geometries \cite{LuTian2004,EnglisZhang2010,LoiMossaZuddas2017,EbenfeltXiaoXu2026}. We are interested in the restrictions that its vanishing imposes on the underlying K\"ahler metric.

To describe this setting, let $M$ be a compact complex manifold of complex dimension $n$, and let $L\rightarrow M$ be a positive holomorphic line bundle equipped with a smooth Hermitian metric $h$ of positive curvature. Set
\[
\omega=\Ric(L,h)=-\ii\pa\pp\log h,
\]
and let $g$ be the associated K\"ahler metric. $(M,\omega;L,h)$ is called a polarization. Denote by $(L^*,h^*)$ the dual Hermitian line bundle and define its unit disk and circle bundles by
\[
D=D(L^*)\coloneq\left\{v\in L^*:\|v\|_{h^*}<1\right\},\qquad S=S(L^*)\coloneq \left\{v\in L^*:\|v\|_{h^*}=1\right\}.
\]
The positivity of the curvature of $h$ implies that $D$ has strictly pseudoconvex boundary $S$, which carries its induced CR structure.

For a smoothly bounded strictly pseudoconvex domain $\Omega$ in $\mathbb{C}^n$, Fefferman's classical asymptotic theory \cite{Fefferman1974} gives a boundary expansion of its Bergman kernel
\[
K_{\Omega}=\frac{\phi}{\rho^{n+1}}+\psi\log \rho,
\]
where $\phi$ and $\psi$ are smooth functions on $\overline{\Omega}$ and $\rho$ is a smooth defining function such that $\Omega=\{\rho>0\}$. The coefficient $\psi$ is usually known as the logarithmic term of the Bergman kernel. When it vanishes to infinite order along the boundary, we say that $\pa \Omega$ is Bergman logarithmically flat.

The notion of Bergman logarithmic flatness is local in nature. Indeed, it depends only on the local CR geometry of the boundary. Consequently, the definition extends to smoothly bounded strictly pseudoconvex domains in arbitrary complex manifolds, and in particular to the unit disk bundle $D\subset L^*$.

To relate the Bergman logarithmic flatness of the unit disk bundle $D$ to the geometry of the polarized manifold, we shall also recall the Bergman kernel function associated with twisted powers of the polarization. Let $(E,h_E)$ be a Hermitian holomorphic line bundle over $M$. For each $k\geq 1$, let $\mathcal{B}_k^E$ denote the Bergman kernel function associated with $H^0(M,L^k\otimes E)$ endowed with the natural $L^2$-inner product induced by $h^k\otimes h_E$ and $\omega$. If $\{s^{(k)}_{\alpha}\}$ is an orthonormal basis, then
\[
\mathcal{B}^E_k(x)=\displaystyle\sum_{\alpha}\|s_{\alpha}^{(k)}(x)\|^2_{h^k\otimes h_E}.
\]
The twisted Bergman kernel admits the well-known Tian--Yau--Zelditch type asymptotic expansion (see, e.g. \cite[Proposition~5.1]{EbenfeltXiaoXu2026})
\begin{equation}\label{eq:asyexp}
\mathcal{B}_k^E(x)\sim \left(\frac{k}{2\pi}\right)^n\sum_{j=0}^{\infty}\frac{a_j^E(x)}{k^j},\qquad \text{uniformly for } x\in M \text{ as } k\rightarrow+\infty.
\end{equation}

When $E$ is the canonical line bundle $K_M$ of $M$ equipped with the standard Hermitian metric induced by $g$, we call $\mathcal{B}_k^{K_M}$ the canonically twisted Bergman kernel. Very recently, Ebenfelt--Xiao--Xu \cite{EbenfeltXiaoXu2026} established a connection between the two types of Bergman kernels introduced above, completing a characterization of the Bergman logarithmic flatness of $D$ in terms of the coefficients in the asymptotic expansion of the canonically twisted Bergman kernel on $M$. Their result can be stated as follows.

\begin{theorem}[{\cite[Theorem~1.12]{EbenfeltXiaoXu2026}}]\label{thm:logflat}
Let $(M,g; L,h)$ be a compact polarized K\"ahler manifold of complex dimension $n$, and let $S$ be the circle bundle of the dual bundle $(L^*,h^*)$. Let $\{a^{K_M}_j\}_{j=0}^{\infty}$ be the coefficient functions appearing in the expansion \eqref{eq:asyexp} of the Bergman kernel function $\mathcal{B}^{K_M}_k$ for $E=K_M$. Then the following are equivalent:

(i) The unit circle bundle $S$ is Bergman logarithmically flat;

(ii) For any $j\geq 0$, $a^{K_M}_{n+2+j}=0$.
\end{theorem}

The above characterization naturally raises the question of rigidity of Bergman logarithmic flatness under deformations of the underlying K\"ahler metric in the same polarization class. More precisely, let $(M,\omega_0; L,h_0)$ be a compact polarized K\"ahler manifold such that the unit circle bundle $S_0$ is Bergman logarithmically flat. We consider the modified Hermitian metric
\[
h_\varphi=h_0e^{-\varphi}, \qquad \omega_\varphi=\omega_0+\ii\pa\pp\varphi>0,\qquad\text{where } \varphi\in C^{\infty}(M)
\]
and ask whether Bergman logarithmic flatness of the corresponding unit circle bundle $S_{\varphi}$ implies that $h_{\varphi}$ is equivalent to $h_0$ up to a holomorphic automorphism and constant rescaling.

In the Szeg\H{o} setting, Lu--Tian \cite[Theorem~1.2]{LuTian2004} established local rigidity on $\mathbb{CP}^n$. More recently, Loi--Palmieri--Zuddas \cite[Theorem~1.1]{LoiPalmieriZuddas2026} extended this result to compact Hermitian symmetric spaces.

\begin{remark}
The local rigidity result of Loi--Palmieri--Zuddas \cite{LoiPalmieriZuddas2026} concerns Szeg\H{o} logarithmic flatness, whereas ours concerns Bergman logarithmic flatness. The associated coefficient equations are different: the former involves the ordinary Tian--Yau--Zelditch coefficients, whereas the latter involves the canonically twisted Bergman coefficients. In particular, their linearizations have different kernels. In our setting, the kernel of the linearization may contain directions beyond those induced by holomorphic automorphisms. We therefore need a nonlinear obstruction argument to exclude these additional directions. 

Bergman logarithmic flatness does not imply Szeg\H{o} logarithmic flatness in general. In Appendix~\ref{app:blowup}, we show that, for every $n\geq2$, the unit circle bundle associated with the product Fubini--Study polarization on $\operatorname{Bl}_p\mathbb{CP}^n$ is Bergman logarithmically flat but not Szeg\H{o} logarithmically flat with respect to the natural measure.
\end{remark}

\subsection{The main results}
Engli\v{s}--Zhang \cite[Theorems~3.3 and~3.5]{EnglisZhang2010} proved that the unit circle bundles associated with invariant positive Hermitian line bundles over irreducible Hermitian symmetric spaces of compact type are Bergman logarithmically flat. More generally, this holds for any polarized compact locally homogeneous K\"ahler manifolds by \cite[Corollary~1.14]{EbenfeltXiaoXu2026}.

This raises the question of whether an analogous local rigidity theorem holds for Bergman logarithmic flatness on compact Hermitian symmetric spaces. In this paper, we answer this question affirmatively, establishing local rigidity near the standard metrics associated with fixed positive polarizations.

We first consider the irreducible case. Let $M$ be an irreducible Hermitian symmetric space of compact type of complex dimension $n$. Its Picard group is generated by the ample line bundle $L=\mathcal O_M(1)$. Let $p$ denote its genus (equivalently its Fano index), and choose an invariant Hermitian metric $h_0$ on $L$. Set
\[
\omega_0=\Ric(L,h_0).
\]
Then $\omega_0$ is the invariant K\"ahler metric in $2\pi c_1(L)$ and satisfies
\[
\Ric(\omega_0)=p\omega_0, \qquad K_M^{-1}\simeq L^p.
\]
With this notation, the main result of this paper is the following local rigidity theorem.

\begin{theorem}[= Theorem~\ref{thm:odd-a-rigidity} for odd $a$; = Theorem~\ref{thm:even-a-rigidity} for even $a$]\label{thm:localrigidity}
Let $M$ be an irreducible Hermitian symmetric space of compact type of complex dimension $n$, with the standard polarization $(M,\omega_0;L,h_0)$ described above. For every $0<\alpha<1$ there are integers $J\geq n+2$ and $s\geq 2$, and $\varepsilon>0$, such that the following holds. If $\varphi\in C^{\infty}(M,\mathbb{R})$,
\[
\omega_{\varphi}>0,\qquad \|\varphi\|_{C^{s,\alpha}(M,\omega_0)}<\varepsilon,
\]
and
\[
a^{K_M}_{n+2}(h_{\varphi})=\cdots=a^{K_M}_{J}(h_{\varphi})=0,
\]
then
\[
h_{\varphi}=CF^*h_0
\]
for some $C>0$ and some polarized automorphism $F\in \operatorname{Aut}^0(M,L)$. In particular, the conclusion holds when $S_{\varphi}$ is Bergman logarithmically flat.
\end{theorem}
A polarized automorphism is a holomorphic automorphism of $M$ that lifts to $L$. In our setting, every automorphism in $\operatorname{Aut}^0(M)$ admits such a lift, and any two lifts of the same automorphism differ by a nonzero constant.

Write $M=G^*/K$, where $G^*$ is a simply connected compact simple Lie group acting transitively on $M$, and $K$ is the isotropy subgroup at a fixed base point. Let $(r,a,b)$ be the restricted-root parameters, with $r$ the rank of $M$. The parameters and the complete list of spaces will be given in Section~\ref{subsec:hermitian-symmetric-spaces}. The $G^*$-action gives the multiplicity-free decomposition
\[
L^2(M,\omega_0)=\widehat{\bigoplus}_{\mathbf s}\mathcal V_{\mathbf s}.
\]
Each $\mathcal V_{\mathbf s}$ is a finite-dimensional irreducible spherical $G^*$-module of smooth complex-valued functions, and each type occurs once. The label $\mathbf s=(s_1,\ldots,s_r)$ is an integer partition with $s_1\geq \cdots\geq s_r\geq 0$. For simplicity, we denote $(s,0,\ldots,0)$ by $(s)$ and $(1,\ldots,1,0,\ldots,0)$, with $q$ ones, by $(1^q)$. Complex conjugation preserves these spaces, so we denote by $\mathcal V_{\mathbf s}^{\mathbb R}$ the subspace of real-valued functions in $\mathcal V_{\mathbf s}$. For further details on this decomposition, see Subsection~\ref{subsec:peter-weyl}.

To understand the local rigidity mechanism behind Theorem~\ref{thm:localrigidity}, we first consider its infinitesimal counterpart. By the characterization Theorem~\ref{thm:logflat} of Ebenfelt--Xiao--Xu, Bergman logarithmic flatness is expressed by the infinite system
\begin{equation}\label{eq: logflatequation}
a^{K_M}_j(h_\varphi)=0,\qquad j\geq n+2.
\end{equation}
Linearizing this system at the standard metric $h_0$ leads to the following infinitesimal rigidity statement.

\begin{theorem}[= Theorem~\ref{thm:linearization-repeat}]\label{thm:linearization}
At the invariant metric $h_0$, the common real kernel of the linearizations of all equations in \eqref{eq: logflatequation} is

(i) $\mathbb{R}\mathbf{1}\oplus \mathcal{V}_{(1)}^{\mathbb{R}}$ when $a$ is odd;

(ii) $\mathbb{R}\mathbf{1}\oplus \mathcal{V}_{(1)}^{\mathbb{R}}\oplus \bigoplus_{q=2}^r \mathcal{V}_{(1^q)}^{\mathbb{R}}$ when $a$ is even.

The first two summands are precisely the infinitesimal constant rescalings and the normalized real holomorphy potentials.
\end{theorem}

As we have seen, the principal difficulty is that the linearized equations may admit infinitesimal deformations which are not induced by constant rescalings or polarized automorphisms. When the common kernel is $\mathbb R\mathbf 1\oplus\mathcal V_{(1)}^{\mathbb R}$, we can establish the local rigidity by using elliptic reduction and the implicit function theorem. When the additional summands $\bigoplus_{q=2}^r\mathcal V_{(1^q)}^{\mathbb R}$ occur, we must show that the corresponding infinitesimal deformations do not extend to actual solutions. The main part of this paper is devoted to establishing higher-order obstructions that exclude these additional infinitesimal directions.

\subsection{Extensions to later coefficient blocks and product spaces}

We consider two extensions of the preceding rigidity result. First, the above rigidity result remains valid if, for any prescribed integer $q\geq n+2$, a sufficiently long finite block of coefficients $a^{K_M}_q,\ldots,a^{K_M}_J$ vanishes.

\begin{theorem}
\label{thm:finite-block-rigidity}
Let $M$ be an irreducible Hermitian symmetric space of compact type of complex dimension $n$, with the standard polarization $(M,\omega_0;L,h_0)$ described above. For every integer $q\geq n+2$ and every $\alpha\in(0,1)$,
there exist integers $J\geq q$ and $s\geq 2$, and a constant $\varepsilon>0$, depending only on $(M,L,h_0)$, $q$, and $\alpha$, such that the following holds. If $\varphi\in C^\infty(M,\mathbb R)$ satisfies
\[
\omega_\varphi>0,\qquad\|\varphi\|_{C^{s,\alpha}(M,\omega_0)}<\varepsilon,
\]
and
\[
a^{K_M}_{q}(h_{\varphi})=\cdots=a^{K_M}_{J}(h_{\varphi})=0,
\]
then
\[
h_\varphi=C\,F^*h_0
\]
for some constant $C>0$ and some $F\in\operatorname{Aut}^0(M,L)$. 
\end{theorem}

In fact, the kernel of the linearized operator remains unchanged when the vanishing conditions start at any index $q\geq n+2$.

Second, we extend the rigidity result to arbitrary compact Hermitian symmetric spaces equipped with a fixed invariant positive polarization. Let $M$ be a connected compact Hermitian symmetric space of complex dimension $n$. Write
\[
M=T\times M_1\times\cdots\times M_l,
\]
where $T$ is a flat complex torus and each $M_i$ is an irreducible Hermitian symmetric space of compact type. For each $i$, let $p_i$ denote the genus of $M_i$, take the ample generator $L_i=\mathcal{O}_{M_i}(1)$, and let $h_i$ be its invariant Hermitian metric, and set $\omega_i=\Ric(L_i,h_i)$. Then
\[
\Ric(\omega_i)=p_i\omega_i,\qquad K_{M_i}^{-1}\simeq L_i^{p_i}.
\]
Every positive line bundle $L$ on $M$ admits a decomposition
\[
L\simeq L_T\boxtimes L_1^{d_1}\boxtimes\cdots\boxtimes L_l^{d_l},\qquad d_i\in\mathbb Z_{>0},
\]
where $L_T$ is a positive line bundle on $T$ and $\boxtimes$ denotes the external tensor product; see Section~\ref{subsec:hermitian-symmetric-spaces}. Choose a Hermitian metric $h_T$ on $L_T$ whose curvature $\omega_T=\Ric(L_T,h_T)$ is translation-invariant, and let
\[
h_0=h_T\boxtimes h_1^{d_1}\boxtimes\cdots\boxtimes h_l^{d_l}.
\]
Let $\pi_T$ and $\pi_i$ be the corresponding factor projections. Then the curvature is
\[
\omega_0=\Ric(L,h_0)=\pi_T^*\omega_T+\sum_{i=1}^ld_i\pi_i^*\omega_i.
\]
Thus $\omega_0$ is a product invariant K\"ahler metric representing $2\pi c_1(L)$.

The following theorem extends the local rigidity Theorem~\ref{thm:finite-block-rigidity} to this more general setting. Note that the deformation potential is arbitrary and need not split across the factors. 

\begin{theorem}\label{thm:generalrigidity}
Let
\[
M=T\times M_1\times\cdots\times M_l
\]
be a connected compact Hermitian symmetric space of complex dimension $n$, equipped with the fixed positive polarization
\[
L=L_T\boxtimes L_1^{d_1}\boxtimes\cdots\boxtimes L_l^{d_l},\qquad h_0=h_T\boxtimes h_1^{d_1}\boxtimes\cdots\boxtimes h_l^{d_l},
\]
as above, where $\Ric(L_T,h_T)$ is translation-invariant and each $h_i$ is invariant. 

For every integer $q\geq n+2$ and every $0<\alpha<1$, there exist integers $J\geq q$ and $s\geq 2$, and $\varepsilon>0$, depending only on $(M,L,h_0)$, $q$, and $\alpha$, such that
the following holds. If $\varphi\in C^\infty(M,\mathbb R)$ satisfies
\[
\omega_\varphi>0,\qquad\|\varphi\|_{C^{s,\alpha}(M,\omega_0)}<\varepsilon,
\]
and
\[
a^{K_M}_q(h_\varphi)=\cdots=a^{K_M}_J(h_\varphi)=0,
\]
then
\[
h_\varphi=CF^*h_0
\]
for some $C>0$ and some polarized automorphism $F\in\operatorname{Aut}^0(M,L)$. In particular, the conclusion holds when $S_\varphi$ is Bergman logarithmically flat. If $M=T$ is a complex torus, one may take $J=q$. Moreover, the conclusion can be strengthened to $h_\varphi=Ch_0$.
\end{theorem}

The proof reduces to the compact-type factors and uses a similar nonlinear obstruction argument. We also determine the common kernel of the linearized operators, which may depend on the polarization and the starting index $q$. See Section~\ref{sec:products} for details.

\subsection{The organization of the proof}
Roughly speaking, our proof begins with an analysis of the linearized Bergman coefficient equations. Since their common kernel may contain infinitesimal directions not induced by holomorphic automorphisms, we establish higher-order obstructions to exclude these additional directions.

Section~\ref{sec:pre} develops the geometric, representation-theoretic, and analytic tools we needed for the proof. The multiplicity-free Peter--Weyl decomposition allows us to describe the Laplacian and the Berezin transform through their eigenvalues on spherical modules. We also find an exact generalized Cartan product identity for Toeplitz operators. This provides the algebraic tool we needed to analyze the nonlinear terms in the Bergman kernel.

In Section~\ref{sec:variation}, we first compute the variation of the Bergman kernel using the Gram matrix, and then obtain the linearization of the corresponding Bergman coefficients. On each spherical summand, these linearizations are expressed in terms of a rational function of $k$, whose Laurent expansion has the form $P(k)+c/k$ when all coefficients of index at least $n+2$ vanish.

In Section~\ref{sec:infinitesimal-deformations}, we use the rational function obtained in Section~\ref{sec:variation} to derive a divisibility condition involving the denominator of the Berezin eigenvalue and the Hilbert polynomial. This yields an explicit classification of the infinitesimal deformations, thereby proving Theorem~\ref{thm:linearization}.

In Section~\ref{sec:elliptic-reduction}, we construct a finite-order elliptic operator whose linearization has the same kernel as the original infinite system. After removing the directions corresponding to constant rescalings and holomorphic automorphisms, the implicit-function theorem gives a finite-dimensional real-analytic graph
\[
\varphi=v+w(v),\qquad v\in E,\qquad w(0)=Dw(0)=0,
\]
where $E$ is the part of the linearized kernel orthogonal to the symmetry directions. Every sufficiently small normalized solution of the original coefficient equations must lie on this graph. 

Section~\ref{sec:odd-a-rigidity} completes the proof of Theorem~\ref{thm:localrigidity} when $a$ is odd. In this case, $E=0$, so the normalized graph consists only of the origin. The same argument also applies in the rank-one case.

Section~\ref{sec:nonlinear-obstructions} establishes the nonlinear obstructions when $a$ is even. Along a nonconstant real-analytic curve $v(t)+w(v(t))$, highest-weight arguments allow us to isolate the leading contribution of a single spherical component. Using the exact Cartan product identity and the inverse Gram matrix recursion, we show that one of the corresponding projected Taylor coefficients of the Bergman kernel is a rational function of $k$ with a nonzero pole at $k=a/2$. This obstruction occurs at the second order when $a\geq4$, while higher orders may be needed when $a=2$.

Section~\ref{sec:local-rigidity-even} proves Theorem~\ref{thm:localrigidity} for even $a$. Using a common denominator and real-analytic curve selection, we show that the origin is isolated among the normalized solutions of finitely many Bergman coefficient equations. The key step is to derive a contradiction with the nonlinear obstruction established in Section~\ref{sec:nonlinear-obstructions}.

In Section~\ref{sec:later-bergman-coefficients}, we extend the argument to coefficient equations starting at any prescribed index $q\geq n+2$ and prove Theorem~\ref{thm:finite-block-rigidity}. Section~\ref{sec:products} extends the results to arbitrary compact Hermitian symmetric spaces with fixed invariant positive polarizations thereby proving Theorem~\ref{thm:generalrigidity}. Finally, we discuss several related problems and further directions in Section~\ref{sec:further-discussion}.

\medskip
\noindent\textbf{Acknowledgements.} The second author is grateful to the Department of Mathematics at the University of California, Irvine, for its hospitality during his visit. After completing our work, we noticed the recent work of Loi--Palmieri--Zuddas \cite{LoiPalmieriZuddas2026}, which establishes a similar local rigidity result for Szeg\H{o} logarithmic flatness on compact Hermitian symmetric spaces. 

\medskip
\noindent\textbf{AI Usage Statement.} ChatGPT (GPT-6 Astra) was used in the preparation of this paper in the following ways:

1. It assisted in gathering and organizing relevant background material on representation theory and Hermitian symmetric spaces, and suggested a proof strategy for Theorem~\ref{thm:pre-exact-cartan}.

2. After the authors established local rigidity of Bergman logarithmic flatness in the irreducible case, they asked ChatGPT to extend their methods to general coefficient blocks and arbitrary compact Hermitian symmetric spaces. ChatGPT successfully developed proofs of these extensions, which were subsequently independently verified and rewritten by the authors.

3. It assisted with language editing and polishing throughout the manuscript.

All mathematical arguments and AI-assisted content were carefully reviewed by the authors, who take full responsibility for the correctness and final content of the paper.

\section{Preliminaries}\label{sec:pre}
This section collects the geometric, representation-theoretic, and analytic tools needed for the rigidity argument. Subsection~\ref{subsec:hermitian-symmetric-spaces} reviews compact Hermitian symmetric spaces and their invariant polarizations, and establishes Bergman logarithmic flatness for the reference metric. Subsection~\ref{subsec:peter-weyl} recalls the Peter--Weyl decomposition into irreducible spherical representations, allowing us to study the linearized equations on each component. Subsection~\ref{subsec:berezin} introduces Toeplitz operators and Berezin transforms and records their spectral formulae, which are used to classify infinitesimal deformations in Section~\ref{sec:infinitesimal-deformations}. Subsection~\ref{subsec:affine-realization} describes the affine realization of compact Hermitian symmetric spaces and their associated Jordan-theoretic structures. We then recall the Fischer inner product, the Hua--Kostant--Schmid decomposition, and the Faraut--Kor\'anyi expansion, culminating in the exact polynomial expansion of Proposition~\ref{prop:pre-compact-fk}. Finally, Subsection~\ref{subsec:pre-conical-cartan} realizes holomorphic sections as polynomials and establishes exact Cartan product identities for Toeplitz operators, providing the key algebraic input for the nonlinear obstruction argument.

\subsection{Hermitian symmetric spaces, restricted-root data, and Bergman logarithmic flatness}\label{subsec:hermitian-symmetric-spaces}

A Hermitian symmetric space is a connected K\"ahler manifold whose geodesic symmetry at each point extends to a global holomorphic isometry. We recall some standard facts about these spaces, following \cite{Helgason1978,Knapp2002}.

We first assume that $M$ is irreducible and of compact type, with complex dimension $n$. It is simply connected, and its Picard group is infinite cyclic, generated by an ample line bundle $L=\mathcal O_M(1)$. Choose an invariant Hermitian metric $h_0$ of positive curvature on $L$, and set $\omega_0=\Ric(L,h_0)$. Write $g_0$ for the associated K\"ahler metric and $J$ for the complex structure. Let $p$ denote the genus of $M$, which coincides with its Fano index. Then
\[
K_M^{-1}\simeq L^p,\qquad \Ric(\omega_0)=p\omega_0.
\]
Let $\operatorname{Aut}^0(M,L)$ denote the identity component of the group of holomorphic automorphisms of $M$ admitting a holomorphic lift to $L$. In the present setting, $\operatorname{Aut}^0(M,L)=\operatorname{Aut}^0(M)$, and any two lifts of the same automorphism to $L$ differ by a nonzero constant.

Let $\mathfrak g$ be a simple Lie algebra with a Cartan decomposition $\mathfrak g=\mathfrak k\oplus\mathfrak p$. Assume it is of Hermitian type; equivalently, the center $\mathfrak z(\mathfrak k)$ is one-dimensional. Choose $Z_0\in\mathfrak z(\mathfrak k)$ such that $\operatorname{ad}(Z_0)^2=-\operatorname{Id}$ on $\mathfrak p$, and let $\mathfrak h\subset\mathfrak k$ be a compact Cartan subalgebra of $\mathfrak g$. Then $Z_0\in\mathfrak h$, and $\mathfrak h^{\mathbb C}$ is a Cartan subalgebra of $\mathfrak g^{\mathbb C}$. Put $H=-iZ_0$. A root $\alpha$ is compact or noncompact according as its root space lies in $\mathfrak k^{\mathbb C}$ or $\mathfrak p^{\mathbb C}$; correspondingly, $\alpha(H)=0$ or $\alpha(H)=\pm1$. Choose a positive root system $\Delta^+$ in which every root satisfying $\alpha(H)=1$ is positive, and let $\Psi^+\subset\Delta^+$ be the set of positive noncompact roots.

Let $\tau$ denote conjugation with respect to the compact real form $\mathfrak g^*=\mathfrak k\oplus i\mathfrak p$. For each root $\alpha$, choose a nonzero root vector $e_\alpha$ such that
\[
\tau(e_\alpha)=-e_{-\alpha},
\qquad
[e_\alpha,e_{-\alpha}]=h_\alpha,
\qquad
[h_\alpha,e_\alpha]=2e_\alpha,
\]
where $h_\alpha\in\mathfrak h^{\mathbb C}$ is the coroot
associated with $\alpha$.
We obtain the Harish--Chandra decomposition
\[
\mathfrak g^{\mathbb C}=\mathfrak p^-\oplus\mathfrak k^{\mathbb C}\oplus\mathfrak p^+,\qquad\mathfrak p^\pm=\bigoplus_{\alpha\in\Psi^+}(\mathfrak g^{\mathbb C})_{\pm\alpha}.
\]
Here $(\mathfrak g^{\mathbb C})_\alpha$ denotes the ordinary root space of $\alpha$. In particular, $\mathfrak p^+$ and $\mathfrak p^-$ are abelian.

Two roots $\alpha,\beta$ with $\alpha\ne\pm\beta$ are called strongly orthogonal if neither $\alpha+\beta$ nor $\alpha-\beta$ is a root. Construct an ordered family $\gamma_1,\ldots,\gamma_r\in\Psi^+$ inductively: let $\gamma_1$ be the largest root in $\Psi^+$, and let $\gamma_j$ be the largest root in $\Psi^+$ orthogonal to $\gamma_1,\ldots,\gamma_{j-1}$, continuing until no such root remains. This is Harish--Chandra's construction of a maximal strongly orthogonal family, and $\mathfrak a=\bigoplus_{j=1}^r\mathbb R(e_{\gamma_j}+e_{-\gamma_j})\subset\mathfrak p$ is a maximal abelian subspace of $\mathfrak p$, with $r=\dim_{\mathbb R}\mathfrak a$ equal to the rank of $M$.

For $\beta\in\mathfrak a^*$, define $\mathfrak g_\beta=\{X\in\mathfrak g\colon [A,X]=\beta(A)X\text{ for every }A\in\mathfrak a\}$. The restricted root system is $\Sigma=\{\beta\in\mathfrak a^*\setminus\{0\}\colon\mathfrak g_\beta\neq 0\}$, and the restricted root decomposition is $\mathfrak g=\mathfrak z_{\mathfrak g}(\mathfrak a)\oplus\bigoplus_{\beta\in\Sigma}\mathfrak g_\beta$, where $\mathfrak z_{\mathfrak g}(\mathfrak a)$ denotes the centralizer of $\mathfrak a$ in $\mathfrak g$. The multiplicity of $\beta$ is
$m_\beta=\dim_{\mathbb R}\mathfrak g_\beta$.

Put $X_j=e_{\gamma_j}+e_{-\gamma_j}$, define the Cayley transform by
\[
C=\exp\left(\frac{\pi}{4}\operatorname{ad}\sum_{j=1}^r(e_{-\gamma_j}-e_{\gamma_j})\right)\in\operatorname{Aut}(\mathfrak g^{\mathbb C}).
\]
The $\mathfrak{sl}_2$ relations and strong orthogonality give $C(h_{\gamma_j})=X_j$. Define $\varepsilon_j=\frac12(\gamma_j\circ C^{-1})\big|_{\mathfrak a}\in\mathfrak a^*$. Then $\varepsilon_i(X_j)=\delta_{ij}$, so $\varepsilon_1,\ldots,\varepsilon_r$ form the basis dual to $X_1,\ldots,X_r$. With this ordering, the positive restricted roots are
\[
\{\varepsilon_i-\varepsilon_j,\ \varepsilon_i+\varepsilon_j\colon1\leq i<j\leq r\}\cup\{2\varepsilon_i\colon1\leq i\leq r\},
\]
possibly together with $\{\varepsilon_i\colon1\leq i\leq r\}$. Their multiplicities are independent of the indices and satisfy
\[
m_{\varepsilon_i\pm\varepsilon_j}=a,\qquad m_{2\varepsilon_i}=1,\qquad m_{\varepsilon_i}=2b.
\]
The complex dimension $n$ and the genus $p$ satisfy
\[
n=r+\frac{a}{2}r(r-1)+br,\qquad p=2+a(r-1)+b.
\]

The irreducible Hermitian symmetric spaces of compact type fall into four classical families and two exceptional types; see \cite[Chapter~X, \S6.3]{Helgason1978}. Their dimensions and parameters are listed in Table~\ref{tab:chss}; see also \cite[Example~3.3]{Ding2022}.

\begin{table}[H]
\centering
\renewcommand{\arraystretch}{1.3}
\begin{tabular}{c|l|ccccc}
\hline
Type & $M$ & $n=\dim_{\mathbb C}M$ & $r$ & $a$ & $b$ & $p$ \\
\hline
AIII & $\operatorname{Gr}(u,u+v)$, $1\le u\le v$
& $uv$ & $u$ & $2$ & $v-u$ & $u+v$ \\
\hline
CI & $\operatorname{Sp}(t)/\operatorname{U}(t)$, $t\ge2$
& $\dfrac{t(t+1)}2$ & $t$ & $1$ & $0$ & $t+1$ \\
\hline
\multirow{2}{*}{DIII}
& $\operatorname{SO}(4t)/\operatorname{U}(2t)$, $t\ge2$
& $t(2t-1)$ & $t$ & $4$ & $0$ & $4t-2$ \\
& $\operatorname{SO}(4t+2)/\operatorname{U}(2t+1)$, $t\ge2$
& $t(2t+1)$ & $t$ & $4$ & $2$ & $4t$ \\
\hline
BDI & $Q^n$, $n\ge3$
& $n$ & $2$ & $n-2$ & $0$ & $n$ \\
\hline
EIII & $E_6/(\operatorname{Spin}(10)\cdot\operatorname{U}(1))$
& $16$ & $2$ & $6$ & $4$ & $12$ \\
\hline
EVII & $E_7/(E_6\cdot\operatorname{U}(1))$
& $27$ & $3$ & $8$ & $0$ & $18$ \\
\hline
\end{tabular}
\caption{Irreducible Hermitian symmetric spaces of compact type and their parameters.}
\label{tab:chss}
\end{table}

More generally, every connected compact Hermitian symmetric space $M$ admits a holomorphic and isometric decomposition
\[
M\simeq T\times M_1\times\cdots\times M_l,
\]
where $T$ is a flat complex torus and $M_1,\ldots,M_l$ are irreducible Hermitian symmetric spaces of compact type. Denote the corresponding projections by $\pi_T$ and $\pi_i$, and let $L_i=\mathcal O_{M_i}(1)$. Since $H^1(M_i,\mathcal O_{M_i})=0$, repeated application of \cite[Chapter~III, Exercise~12.6]{Hartshorne1977} yields $\operatorname{Pic}(M)\simeq\operatorname{Pic}(T)\oplus\bigoplus_{i=1}^{l}\mathbb Z[L_i]$. Consequently, every positive holomorphic line bundle $L$ on $M$ has the form
\begin{equation}\label{eq:general-product-polarization}
L\simeq L_T\boxtimes L_1^{d_1}\boxtimes\cdots\boxtimes L_l^{d_l},\qquad d_i\in\mathbb Z_{>0},
\end{equation}
with $L_T$ positive on $T$. Here $\boxtimes$ denotes the external tensor product.

Choose invariant Hermitian metrics $h_i$ on $L_i$ and put $\omega_i=\Ric(L_i,h_i)$. On $T$, let $\omega_T$ be the unique translation-invariant representative of $2\pi c_1(L_T)$. The K\"ahler metric defined by $\omega_T$ is flat. Set
\begin{equation}\label{eq:general-reference-metric}
\begin{split}
h_0&=h_T\boxtimes h_1^{d_1}\boxtimes\cdots\boxtimes h_l^{d_l},\\
\omega_0&=\Ric(L,h_0)=\pi_T^*\omega_T+\sum_{i=1}^{l}d_i\pi_i^*\omega_i.
\end{split}
\end{equation}
Thus $\omega_0$ is a symmetric K\"ahler form in $2\pi c_1(L)$. If $p_i$ is the Fano index of $M_i$, then
\[
K_M^{-1}\simeq\mathcal O_T\boxtimes L_1^{p_1}\boxtimes\cdots\boxtimes L_l^{p_l}, \qquad \Ric(\omega_0)=\sum_{i=1}^{l}p_i\pi_i^*\omega_i.
\]

We use $h_0$ as the reference metric for this polarization. A general smooth polarized Hermitian metric on $L$ is written as
\[
h_\varphi=h_0e^{-\varphi},\qquad \omega_\varphi=\omega_0+\ii\pa\pp\varphi>0, \qquad \varphi\in C^\infty(M,\mathbb R).
\]
The potential $\varphi$ may depend jointly on all factors.

The reference metrics above are Bergman logarithmically flat by \cite{EnglisZhang2010,EbenfeltXiaoXu2026}. For completeness, we recall their argument.

\begin{proposition}[see {\cite[Theorem~3.5]{EnglisZhang2010}} and {\cite[Corollary~1.14]{EbenfeltXiaoXu2026}}]
\label{prop:symmetric-reference-log-flat}
For the symmetric reference metric $h_0$, the canonically twisted Bergman coefficients satisfy
\[
a^{K_M}_j(h_0)\equiv0\qquad(j\ge n+1).
\]
Consequently, the unit circle bundle of $(L^*,h_0^*)$ is Bergman logarithmically flat.
\end{proposition}

\begin{proof}
Consider the Bergman kernel $\mathcal B^{K_M}_k(h_0;\cdot)$ associated with $H^0(M,L^k\otimes K_M)$, where $K_M$ carries the metric induced by $\omega_0$ and the $L^2$ inner product is defined using the normalized volume form $d\mu_0=\frac{\omega_0^n}{\int_M\omega_0^n}$. Write its expansion as 
\[
\mathcal B^{K_M}_k(h_0;x)\sim\left(\frac{k}{2\pi}\right)^n\sum_{j=0}^{\infty}a^{K_M}_j(h_0;x)k^{-j},
\]
here the coefficients $a^{K_M}_j(h_0;\cdot)$ are universal polynomials in the curvature tensor and its covariant derivatives (see \cite[Theorem~1.1]{Lu2000} and \cite[Proposition~5.1]{EbenfeltXiaoXu2026}). Since $(M,\omega_0)$ is homogeneous, these coefficients are constant on $M$. Kodaira vanishing gives us
\[
\int_M \mathcal B^{K_M}_k(h_0;x)\,d\mu_0(x)=\operatorname{dim} H^0(M,L^k\otimes K_M)=\chi(M,L^k\otimes K_M).
\]
By the Hirzebruch--Riemann--Roch theorem, this holomorphic Euler characteristic is a polynomial in $k$ of degree $n$. Integrating the asymptotic expansion and comparing coefficients then yields
\[
\int_M a^{K_M}_j(h_0;x)\,d\mu_0(x)=0\qquad(j>n).
\]
Since these coefficients are constant, they vanish identically. In particular, $a^{K_M}_j(h_0)=0$ for all $j\geq n+2$.
\end{proof}

\subsection{Spherical representations and Peter--Weyl decomposition}\label{subsec:peter-weyl}

Let $M=G^*/K$ be an irreducible Hermitian symmetric space of compact type, where $G^*$ is the simply connected compact group with Lie algebra $\mathfrak g^*=\mathfrak k\oplus i\mathfrak p$. 

Let $\Pi$ be the set of simple roots corresponding to $\Delta^+$. For a finite-dimensional complex $G^*$-representation $(\pi,V)$, its differential extends complex linearly to $\mathfrak g^{\mathbb C}$. For $\lambda\in(\mathfrak h^{\mathbb C})^*$, define
\[
V[\lambda]=\{v\in V\colon d\pi(H)v=\lambda(H)v\text{ for all }H\in\mathfrak h^{\mathbb C}\}.
\]
We call $\lambda$ a weight of $V$ if $V[\lambda]\neq\{0\}$. A nonzero vector $v_\lambda\in V[\lambda]$ is called a highest-weight vector if $d\pi(e_\alpha)v_\lambda=0$ for all $\alpha\in\Delta^+$. A weight is called dominant integral if $\lambda(h_\alpha)\in\mathbb Z_{\geq 0}$ for every $\alpha\in\Pi$.

Irreducible finite-dimensional complex $G^*$-representations are classified by dominant integral highest weights. We write $V_\lambda$ for the representation of highest weight $\lambda$, whose highest-weight space is one-dimensional and whose other weights are obtained from $\lambda$ by subtracting nonzero nonnegative integral combinations of simple roots; see \cite[Theorems~5.5 and~5.110]{Knapp2002}. An irreducible $G^*$-representation $(\pi,V)$ is called $K$-spherical if
\[
V^K=\{v\in V\colon\pi(k)v=v\text{ for every }k\in K\}\neq\{0\}.
\]
For a symmetric pair, this fixed-vector space has dimension at most one. We describe the spherical highest weights using the strongly orthogonal roots $\gamma_1,\ldots,\gamma_r$ and the Cayley transform $C$ introduced above. Set $\mathfrak t=\{H\in\mathfrak h\colon\gamma_j(H)=0,\ 1\leq j\leq r\}$. Then $C(\mathfrak h^{\mathbb C})=\mathfrak t^{\mathbb C}\oplus\mathfrak a^{\mathbb C}$ and $C|_{\mathfrak t^{\mathbb C}}=I$. For a highest weight $\lambda$, write
\[
\widetilde\lambda=\lambda\circ C^{-1},\qquad\lambda_{\mathfrak a}=\widetilde\lambda|_{\mathfrak a}.
\]
By the Cartan--Helgason theorem \cite[Theorem~8.49]{Knapp2002}, $V_\lambda$ is $K$-spherical if and only if
\[
\widetilde\lambda|_{\mathfrak t}=0,\qquad\frac{(\lambda_{\mathfrak a},\beta)}{(\beta,\beta)}\in\mathbb Z\quad\text{for every }\beta\in\Sigma,
\]
where the inner product on $\mathfrak a^*$ is induced by the Killing form.

In our coordinates, the spherical highest weights are therefore exactly
\[
\Lambda_{\mathbf s}=\sum_{j=1}^r s_j\gamma_j,\qquad \mathbf s=(s_1,\ldots,s_r)\in\mathcal P_r,
\]
where
\[
\mathcal P_r=\{(s_1,\ldots,s_r)\in\mathbb Z^r\colon s_1\geq\cdots\geq s_r\geq 0\}.
\]

We abbreviate $(s,0,\ldots,0)$ by $(s)$ and the partition with $q$ ones and $r-q$ zeros by $(1^q)$. The Peter--Weyl theorem decomposes $L^2(G^*)$ into matrix coefficients of irreducible representations; see \cite[Theorem~4.20]{Knapp2002}. On the quotient it yields the following decomposition.

\begin{theorem}[see {\cite[Proposition~3.1]{Zhang1998}}]
\label{thm:spherical-decomposition}
Under the $G^*$-action, the space $L^2(M,\mathbb C)$ admits the orthogonal decomposition
\[
L^2(M,\mathbb C)=\widehat{\bigoplus}_{\mathbf s\in\mathcal P_r}\mathcal V_{\mathbf s}.
\]
Here $\mathcal V_{\mathbf s}\subset C^\infty(M,\mathbb C)$ is an irreducible $K$-spherical $G^*$-module of highest weight $\Lambda_{\mathbf s}=\sum_{j=1}^r s_j\gamma_j$, and each representation occurs with multiplicity one. Moreover, $\operatorname{dim}_{\mathbb C}\mathcal V_{\mathbf s}^{K}=1$.
\end{theorem}
Each $\mathcal V_{\mathbf s}$ is preserved by complex conjugation. Set $\mathcal V_{\mathbf s}^{\mathbb R}=\mathcal V_{\mathbf s}\cap C^\infty(M,\mathbb R)$. Then $\mathcal V_{\mathbf s}=\mathcal V_{\mathbf s}^{\mathbb R}\otimes_{\mathbb R}\mathbb C$.

Recall that a real-valued smooth function $f$ is called a real holomorphy potential if $\nabla^{1,0}f$ is holomorphic, that is, $\pp\nabla^{1,0}f=0$.

\begin{proposition}\label{prop:constants-holomorphy}
\[
\mathcal V_{\mathbf0}^{\mathbb R}=\mathbb R\mathbf1,\qquad
\mathcal V_{(1)}^{\mathbb R}=\left\{f\in C^\infty(M,\mathbb R)\colon \pp\nabla^{1,0}f=0,\quad \int_M f\,d\mu_0=0\right\}.
\]
In particular, the full space of real holomorphy potentials is $\mathcal S=\mathbb R\mathbf1\oplus\mathcal V_{(1)}^{\mathbb R}$.
\end{proposition}
\begin{proof}
The highest weight $\Lambda_{\mathbf0}=0$ corresponds to the trivial representation. Its functions are $G^*$-invariant and therefore constant, since $G^*$ acts transitively on $M$.

For $X\in\mathfrak g^*=\operatorname{Lie} (G^*)$, let $X_M$ be the vector field generated by the one-parameter group $\exp(tX)$. It is a real holomorphic Killing vector field. Since $M$ is simply connected, there exists a unique real function $f_X$ satisfying
\[
df_X=-\iota_{X_M}\omega_0,\qquad\int_M f_X\,d\mu_0=0,
\]
which gives $X_M=J\nabla f_X$. Since $X_M$ is holomorphic, $f_X$ is a real holomorphy potential.

Conversely, if $f$ is a real holomorphy potential, then $J\nabla f$ is holomorphic and $\iota_{J\nabla f}\omega_0=-df$. It is therefore a holomorphic Killing field and belongs to $\mathfrak{g}^*$. After normalization, every real holomorphy potential is consequently of the form $f_X$.

The correspondence $X\mapsto f_X$ is $G^*$-equivariant because
\[
f_{\operatorname{Ad}(g)X}=g\cdot f_X.
\]
Its complexification therefore identifies the complexified space of normalized real holomorphy potentials with the adjoint representation $(\mathfrak g^*)^{\mathbb C}=\mathfrak g^{\mathbb C}$. It is irreducible with highest weight $\gamma_1=\Lambda_{(1)}$. Therefore the proof is complete by taking the real form.
\end{proof}

\subsection{Toeplitz operators and Berezin transform}
\label{subsec:berezin}

Toeplitz operators and Berezin transforms play a central role in the quantization of K\"ahler manifolds; see \cite{Berezin1974,BordemannMeinrenkenSchlichenmaier1994,KarabegovSchlichenmaier2001,MaMarinescu2012} and the survey \cite{Schlichenmaier2010}. We recall their definitions and describe their action on the spherical modules $\mathcal V_{\mathbf s}$ here. We work with the standard polarization $(M,L,h_0)$ and the invariant probability measure $d\mu_0$ introduced above. For an integer $m\geq 0$, equip $\mathcal H_m=H^0(M,L^m)$ with the Hermitian inner product
\[
\langle u,v\rangle_m=\int_M h_0^m(u,v)\,d\mu_0,
\]
let
\[
P_m\colon L^2(M,L^m)\longrightarrow\mathcal H_m
\]
be the orthogonal projection. For $f\in C^\infty(M,\mathbb C)$, its Toeplitz operator is
\[
T_m(f)u=P_m(fu),\qquad \text{for }u\in\mathcal H_m.
\]

Let $K_m(x,y)$ denote the reproducing kernel of $\mathcal H_m$, viewed as a section of $L^m\boxtimes(L^m)^*$ over $M\times M$. For $A\in\operatorname{End}(\mathcal H_m)$, define its covariant symbol by
\[
\sigma_m(A)(x)=\frac{(AK_m(\cdot,x))(x)}{K_m(x,x)}.
\]
The quotient is well-defined since $L_x^m$ is one-dimensional and $K_m(x,x)\neq 0$. The Berezin transform is defined by
\[
B_m f=\sigma_m(T_m(f)).
\]
These constructions are reviewed in \cite[Sections~3, 6, and~7]{Schlichenmaier2010}. Writing $K_m(x,x)=\rho_m(x)\operatorname{Id}_{L_x^m}$, the Berezin transform admits the integral representation
\[
B_m f(x)=\frac{1}{\rho_m(x)}\int_M \|K_m(x,y)\|^2 f(y)\,d\mu_0(y),
\]
where the norm is induced by $h_0^m$ on the two fibers.

Homogeneity implies that $\rho_m$ is constant. Since $d\mu_0$ has total mass one, we obtain $\rho_m=\dim_{\mathbb C}\mathcal H_m$.
By definition, $B_m\mathbf{1}=\mathbf{1}$. The invariance of the metric and measure also implies that $B_m$ commutes with the $G^*$-action.

To express the eigenvalues of $B_m$, set
\[
c_j=\frac a2(j-1),\qquad 1\leq j\leq r.
\]
For $\mathbf s\in\mathcal P_r$, put $|\mathbf s|=\sum_{j=1}^r s_j$ and define the generalized Pochhammer symbol by
\[
(z)_{\mathbf s}=\prod_{j=1}^r(z-c_j)_{s_j},
\]
where $(u)_k=\prod_{l=0}^{k-1}(u+l)$ for $k\geq 1$, with $(u)_0=1$.

Since $B_m$ is $G^*$-equivariant and the spherical decomposition is multiplicity-free, Schur's lemma then implies that $B_m$ acts by a scalar on each irreducible summand $\mathcal V_{\mathbf s}$. The following formula is due to Zhang \cite{Zhang1998}.

\begin{theorem}[{\cite[Theorem~3.6]{Zhang1998}}]
\label{thm:berezin-eigenvalues}
For $\mathbf s\in\mathcal P_r$, define the polynomials
\begin{equation}\label{eq:numerator}
A_{\mathbf s}(m)=(-1)^{|\mathbf s|}(-m)_{\mathbf s}=\prod_{j=1}^r\prod_{l=0}^{s_j-1}(m+c_j-l),
\end{equation}
and
\begin{equation}\label{eq:denominator}
D_{\mathbf s}(m)=(m+p)_{\mathbf s}=\prod_{j=1}^r\prod_{l=0}^{s_j-1}(m+p-c_j+l),\qquad c_j=\frac{a}{2}(j-1),
\end{equation}
Both polynomials are monic of degree $|\mathbf s|$.

Then for every integer $m\ge0$, the Berezin transform acts on $\mathcal V_{\mathbf s}$ by
\begin{equation}\label{eq:berezin-eigenvalues}
B_m|_{\mathcal V_{\mathbf s}}=b_{\mathbf s}(m)\,I,\qquad \text{where }b_{\mathbf s}(m)=\frac{A_{\mathbf s}(m)}{D_{\mathbf s}(m)}.
\end{equation}
Moreover, for any nonnegative integer $m$, $D_{\mathbf s}(m)>0$ and,
\[
b_{\mathbf s}(m)>0 \quad\text{if }s_1\leq m,\qquad b_{\mathbf s}(m)=0 \quad\text{if }s_1>m.
\]
In particular,
\[
\operatorname{Ran}B_m=\bigoplus_{\substack{\mathbf s\in\mathcal P_r\\s_1\le m}}\mathcal V_{\mathbf s}.
\]
\end{theorem}

\begin{remark}
For each fixed $\mathbf s\in\mathcal P_r$, the right-hand side of \eqref{eq:berezin-eigenvalues} defines a rational function of the complex parameter $m$. This observation will be useful in the sequel, as it allows us to treat the eigenvalue formula as a rational identity in $m$ and to exploit its algebraic properties beyond the original range of positive integer values.
\end{remark}

We next consider the action of the Laplacian. Since $\Delta$ is $G^*$-invariant, Schur's lemma implies that $\Delta$ also acts by a scalar on each $\mathcal V_{\mathbf s}$.

\begin{theorem}\label{thm:laplace-spectrum}
For every $\mathbf s\in\mathcal P_r$,
\[
\Delta|_{\mathcal V_{\mathbf s}}=\lambda_{\mathbf s}I,\qquad\text{where }\lambda_{\mathbf s}=\sum_{j=1}^r s_j\bigl(s_j+p-1-a(j-1)\bigr).
\]
In particular,
\[
\lambda_{\mathbf0}=0,\qquad\lambda_{(1)}=p,\qquad\lambda_{(1^q)}=qp-\frac a2q(q-1)\quad(1\leq q\leq r).
\]
\end{theorem}

\begin{proof}
By \cite{Englis2000} and \cite[Theorem~7.2]{Schlichenmaier2010}, the Berezin transform has the asymptotic expansion
 \[
B_mf=f-\frac1m\Delta f+O(m^{-2});
\]
For a fixed partition $\mathbf s$, each factor in \eqref{eq:berezin-eigenvalues} satisfies
\[
\frac{m+c_j-l}{m+p-c_j+l}=1-\frac{p-2c_j+2l}{m}+O(m^{-2}).
\]
Multiplying these finitely many factors gives
\[
b_{\mathbf s}(m)=1-\frac1m\sum_{j=1}^r\sum_{l=0}^{s_j-1}(p-2c_j+2l)+O(m^{-2}).
\]
The double sum equals
\[
\sum_{j=1}^rs_j\bigl(s_j+p-1-a(j-1)\bigr).
\]
Comparison of the two expansions on $\mathcal V_{\mathbf s}$ proves the assertion.
\end{proof}

\begin{remark}
The first nonzero eigenvalue of $\Delta$ is $p$, and its eigenspace is precisely $\mathcal V_{(1)}$. Indeed, if $\mathbf s\ne\mathbf0$, then $s_1\geq1$, so $s_1(s_1+p-1)\geq p$, with equality only when $s_1=1$. All remaining summands are nonnegative and are strictly
positive whenever $s_j>0$, since $p-1-a(j-1)=1+b+a(r-j)>0$. Thus $\lambda_{\mathbf s}=p$ holds only for $\mathbf s=(1)$.
\end{remark}

In particular, $B_m$ and $\Delta$ commute on $C^\infty(M)$, since they both are scalar on each spherical summand. The operators above are defined using sections of $L^m$. We note that since $K_M\simeq L^{-p}$, $L^k\otimes K_M\simeq L^{k-p}$. At the invariant metric, the induced Hermitian metrics under this identification differ only by a positive constant.

\subsection{The affine realization and the Faraut--Kor\'anyi expansion}
\label{subsec:affine-realization}
We recall the affine realization of compact Hermitian symmetric spaces and the associated Jordan-theoretic structures, following \cite{Knapp2002,Loos1977,Schwarz2013}. The space $Z=\mathfrak p^+$ carries a positive Hermitian Jordan triple product defined by
\[
\{x,y,z\}\coloneq-[[x,\tau y],z],
\]
where $\tau$ is conjugation with respect to the compact real form $g^*$. Write $D(x,y)z=\{x,y,z\}$ and $Q(x)y=\frac12\{x,y,x\}$. The Bergman operator and generic norm are defined by
\[
B(x,y)=I-D(x,y)+Q(x)Q(y),\qquad \det B(z,w)=h(z,w)^p,
\]
with $h(z,0)=h(0,w)=1$. Here $h(z,w)$ is an irreducible polynomial in $z$ and $\overline w$; see \cite[\S4.15]{Loos1977} and \cite[(2.5)--(2.6)]{Schwarz2013}. Equip $Z$ with the normalized trace form
\[
\langle x,y\rangle_0=\frac1p\operatorname{tr}_Z D(x,y),
\]
then on the affine chart, the invariant Hermitian metric and K\"ahler form of $M$ are given by
\[
\|e(z)\|_{h_0}^2=h(z,-z)^{-1},\qquad \omega_0=\ii\pa\pp\log h(z,-z),
\]
where $e$ is a holomorphic frame of $L$, and
\[
d\mu_0(z)=c\,h(z,-z)^{-p}\,dV(z),
\]
where $dV$ is Lebesgue measure on $Z$ and $c>0$; see \cite{EnglisZhang2010}. With respect to orthonormal coordinates for $\langle\cdot,\cdot\rangle_0$, the Fischer inner product on the polynomial algebra $\mathcal P(Z)$ is defined by
\[
\langle z^\alpha,z^\beta\rangle_F=\alpha!\,\delta_{\alpha\beta}.
\]
The Fischer--Riesz map is the conjugate-linear identification
\[
\mathcal R\colon\mathcal P(Z)\longrightarrow\operatorname{Sym}(Z),\qquad ((\mathcal Rq)f)(0)=\langle f,q\rangle_F,
\]
where $\operatorname{Sym}(Z)$ is identified with constant-coefficient holomorphic differential operators. In particular, if $q(z)=\sum_\alpha q_\alpha z^\alpha$, then $\mathcal Rq=\sum_\alpha\overline{q_\alpha}\,\pa^\alpha$.

We recall the Hua--Kostant--Schmid decomposition in the following form.
\begin{theorem}[Hua--Kostant--Schmid decomposition {\cite{Hua1963,Schmid1969,FarautKoranyi1990}}]
\label{thm:pre-hks}
The polynomial algebra $\mathcal P(Z)$ has the orthogonal, multiplicity-free decomposition
\[
\mathcal P(Z)=\bigoplus_{\mathbf s\in\mathcal P_r}\mathcal P_{\mathbf s}(Z).
\]
Each $\mathcal P_{\mathbf s}(Z)$ is an irreducible $K^{\mathbb C}$-module consisting of homogeneous polynomials of total degree $|\mathbf s|$. The corresponding constant-tensor module
\[
\mathcal D_{\mathbf s}\coloneq\mathcal R\mathcal P_{\mathbf s}(Z)\subset\operatorname{Sym}^{|\mathbf s|}(Z)
\]
has highest weight $\Lambda_{\mathbf s}=\sum_j s_j\gamma_j$. Equivalently, $\mathcal P_{\mathbf s}(Z)$ has lowest weight $-\Lambda_{\mathbf s}$.
\end{theorem}

Let $K_{\mathbf s}^F(z,w)$ be the Fischer reproducing kernel of $\mathcal P_{\mathbf s}(Z)$. Orthogonality gives
\begin{equation}\label{eq:pre-fischer-contraction}
\left.(\mathcal Rq)_zK_{\mathbf t}^F(z,w)\right|_{z=0}=
\begin{cases}
\overline{q(w)},&\mathbf t=\mathbf s,\\
0,&\mathbf t\ne\mathbf s,
\end{cases}
\qquad \text{where }q\in\mathcal P_{\mathbf s}(Z).
\end{equation}
Let $\Omega\subset Z$ be the Harish-Chandra realization of the noncompact dual. The Faraut--Kor\'anyi expansion \cite[Theorem~3.8]{FarautKoranyi1990} states that for every $\nu\in \mathbb C$ and every $z,w\in \Omega\subset Z$,
\begin{equation}\label{eq:pre-fk}
h(z,w)^{-\nu}=\sum_{\mathbf s\in\mathcal P_r}(\nu)_{\mathbf s}K_{\mathbf s}^F(z,w).
\end{equation}
Substituting $\nu=-m$ and replacing $w$ by $-w$ in \eqref{eq:pre-fk}, we obtain the following finite expansion on the compact dual.
\begin{proposition}\label{prop:pre-compact-fk}
For every nonnegative integer $m$,
\begin{equation}\label{eq:pre-compact-fk}
h(z,-w)^m=\sum_{s_1\le m}A_{\mathbf s}(m)K_{\mathbf s}^F(z,w),\qquad z,w\in Z,
\end{equation}
where $A_{\mathbf s}(m)=(-1)^{|\mathbf s|}(-m)_{\mathbf s}$ is defined in \eqref{eq:numerator}.
\end{proposition}
\begin{proof}
Set $\nu=-m$ and replace $w$ by $-w$ in \eqref{eq:pre-fk}. Bihomogeneity contributes the factor $(-1)^{|\mathbf s|}$, while $(-m)_{s_1}=0$ for $s_1>m$. Then the resulting finite identity extends from $\Omega\times\Omega$ to $Z\times Z$ as a polynomial identity.
\end{proof}

\subsection{Generalized Cartan products}
\label{subsec:pre-conical-cartan}
In this subsection, we prove an exact Cartan product formula for Toeplitz operators, using the polynomial realization of holomorphic sections. This provides the algebraic foundation for the nonlinear obstruction argument in Section~\ref{sec:nonlinear-obstructions}.

By the Borel--Weil theorem, the space $\mathcal H_m=H^0(M,L^m)$ is an irreducible $G^{\mathbb C}$-representation for every nonnegative integer $m$. Choose $e_o\in L_o$ with $\|e_o\|_{h_0}=1$. Let $\widetilde g$ denote the lifted action of $g\in G^{\mathbb C}$ on $L$. This gives a holomorphic frame $e(z)=\widetilde{\exp(z)}\,e_o$ over the affine chart $Z$. Writing $s(z)=f(z)e(z)^{\otimes m}$, we identify sections with their coefficient functions. Since $Z=\mathfrak p^+$ is abelian, the derived representation of $Z$ acts by
\[
d\pi_m(X)f=-\pa_Xf,\qquad X\in Z.
\]
In the frame $e(z)$, if we write $K_m(z,w)=K_m^{\operatorname{loc}}(z,w)\,e(z)^{\otimes m}\otimes(e(w)^{\otimes m})^*$, then
\[
K_m^{\operatorname{loc}}(z,w)=\rho_m h(z,-w)^m,\qquad \rho_m=\dim_{\mathbb C}\mathcal H_m.
\]
Since $h(z,-w)^m$ is polynomial in $z$, the coefficient functions of all sections in $\mathcal H_m$ are polynomials. Moreover, by Proposition~\ref{prop:pre-compact-fk}, the expansion of $h(z,-w)^m$ involves precisely the polynomial components $\mathcal P_{\mathbf s}(Z)$ with $s_1\leq m$, thus under the identification, we obtain
\[
\mathcal H_m=\bigoplus_{s_1\leq m}\mathcal P_{\mathbf s}(Z).
\]

We now choose highest vectors compatible with Cartan products. For $1\leq j\leq r$, fix a nonzero highest vector $Q_{(1^j)}\in\mathcal D_{(1^j)}$, and set
\[
Q_{\mathbf s}=\prod_{j=1}^r Q_{(1^j)}^{\,s_j-s_{j+1}},\qquad s_{r+1}=0,\qquad q_{\mathbf s}=\mathcal R^{-1}Q_{\mathbf s}.
\]
Here $Q_{\mathbf s}\in\mathcal D_{\mathbf s}$ and $q_{\mathbf s}\in\mathcal P_{\mathbf s}(Z)$ are nonzero highest and lowest weight vectors, respectively. Moreover,
\[
Q_{\mathbf s}Q_{\mathbf t}=Q_{\mathbf s+\mathbf t},\qquad q_{\mathbf s}q_{\mathbf t}=q_{\mathbf s+\mathbf t}.
\]
The derived representation extends to $\operatorname{Sym}(Z)$ via the Poincar\'e--Birkhoff--Witt theorem. We define
\[
C_{\mathbf s}^{(m)}\coloneq(-1)^{|\mathbf s|}d\pi_m(Q_{\mathbf s})=(\mathcal Rq_{\mathbf s})|_{\mathcal H_m}.
\]

Then $C_{\mathbf s}^{(m)}C_{\mathbf t}^{(m)}=C_{\mathbf s+\mathbf t}^{(m)}$. We now relate these operators to their covariant symbols. The covariant symbol map $\sigma_m$ is equivariant for the conjugation action of $G^*$ on $\operatorname{End}(\mathcal H_m)$. Every function $f\in\mathcal V_{\mathbf s}$ is real analytic and admits a unique local polarization $f^{\mathrm{pol}}(z,w)$, holomorphic in $z$ and antiholomorphic in $w$, with $f^{\mathrm{pol}}(z,z)=f(z)$. Define
\[
\mathcal J^-f(w)\coloneq f^{\mathrm{pol}}(0,w).
\]
For $A\in\operatorname{End}(\mathcal H_m)$, the polarized covariant symbol is given near $(0,0)$ by
\[
\sigma_m(A)^{\mathrm{pol}}(z,w)=\frac{(AK_{m,w}^{\operatorname{loc}})(z)}{K_m^{\operatorname{loc}}(z,w)},\qquad \text{where }K_{m,w}^{\operatorname{loc}}=K_m^{\operatorname{loc}}(\cdot,w).
\]
\begin{lemma}\label{lem:pre-symbol-injective}
The covariant symbol map $\sigma_m$ is injective.
\end{lemma}
\begin{proof}
If $\sigma_m(A)=0$, then $(AK_{m,w}^{\operatorname{loc}})(z)$ vanishes on the diagonal and hence vanishes identically. Thus $AK_{m,w}^{\operatorname{loc}}=0$ for every $w$. Since $K_{m,w}^{\operatorname{loc}}$ span the whole $\mathcal H_m$, we obtain $A=0$.
\end{proof}
\begin{proposition}\label{prop:pre-highest-symbols}
There are highest vectors $e_{\mathbf s}\in\mathcal V_{\mathbf s}$, independent of $m$, uniquely specified by
\[
\mathcal J^-e_{\mathbf s}=\overline{q_{\mathbf s}}.
\]
For every integer $m\geq s_1$, the operator $C_{\mathbf s}^{(m)}$ is nonzero and
\[
\sigma_m(C_{\mathbf s}^{(m)})=A_{\mathbf s}(m)e_{\mathbf s}.
\]
Moreover, $e_{\mathbf0}=1$ and $e_{\mathbf s}e_{\mathbf t}=e_{\mathbf s+\mathbf t}$.
\end{proposition}
\begin{proof}
Every positive compact root vector annihilates $Q_{\mathbf s}$, while every positive noncompact root vector commutes with it because $Z$ is abelian. Hence $C_{\mathbf s}^{(m)}$, if nonzero, is a highest vector of weight $\Lambda_{\mathbf s}$ for the conjugation representation on $\operatorname{End}(\mathcal H_m)$. By equivariance, its nonzero symbol belongs to the one-dimensional highest-weight space of $\mathcal V_{\mathbf s}$.

Applying $\mathcal Rq_{\mathbf s}$ in the first variable of the reproducing kernel, we obtain from \eqref{eq:pre-compact-fk} and \eqref{eq:pre-fischer-contraction} that
\[
\begin{aligned}
(\mathcal Rq_{\mathbf s})_z K_m^{\operatorname{loc}}(z,w)\big|_{z=0}&=\rho_m\sum_{t_1\leq m}A_{\mathbf t}(m)
(\mathcal Rq_{\mathbf s})_z K_{\mathbf t}^F(z,w)\big|_{z=0}\\
&=\rho_m A_{\mathbf s}(m)\overline{q_{\mathbf s}(w)}.
\end{aligned}
\]
Since $K_m^{\operatorname{loc}}(0,w)=\rho_m$, the formula for the polarized covariant symbol gives
\[
\mathcal J^-\sigma_m(C_{\mathbf s}^{(m)})(w)
=A_{\mathbf s}(m)\overline{q_{\mathbf s}(w)}.
\]

For $m\geq s_1$, the right-hand side is nonzero. Dividing by $A_{\mathbf s}(m)$ gives a highest vector satisfying the required normalization. Since $\mathcal J^-$ is injective on the one-dimensional highest-weight space, this vector is unique and independent of $m$.
Finally, $e_{\mathbf s}e_{\mathbf t}$ is a highest vector of weight $\Lambda_{\mathbf s+\mathbf t}$, and
\[
\mathcal J^-(e_{\mathbf s}e_{\mathbf t})=\overline{q_{\mathbf s}q_{\mathbf t}}=\overline{q_{\mathbf s+\mathbf t}}.
\]
Uniqueness gives $e_{\mathbf s}e_{\mathbf t}=e_{\mathbf s+\mathbf t}$.
\end{proof}
Let $\operatorname{pr}_{\mathbf s}$ denote the $L^2$-orthogonal projection onto $\mathcal V_{\mathbf s}$. Define the Cartan product by
\[
f\odot g\coloneq\operatorname{pr}_{\mathbf s+\mathbf t}(fg),\qquad \text{for }f\in\mathcal V_{\mathbf s},\quad g\in\mathcal V_{\mathbf t}.
\]
It is nonzero and satisfies $e_{\mathbf s}\odot e_{\mathbf t}=e_{\mathbf s+\mathbf t}$.
\begin{theorem}\label{thm:pre-exact-cartan}
For $f\in\mathcal V_{\mathbf s}$, $g\in\mathcal V_{\mathbf t}$, and every integer $m\geq s_1+t_1$,
\[
\operatorname{pr}_{\mathbf s+\mathbf t}\sigma_m\bigl(T_m(f)T_m(g)\bigr)=\frac{A_{\mathbf s+\mathbf t}(m)}{D_{\mathbf s}(m)D_{\mathbf t}(m)}\,(f\odot g).
\]
\end{theorem}
\begin{proof}
Both sides define $G^*$-equivariant bilinear maps from $\mathcal V_{\mathbf s}\times\mathcal V_{\mathbf t}$ to $\mathcal V_{\mathbf s+\mathbf t}$. Since the Cartan component occurs with multiplicity one in $\mathcal V_{\mathbf s}\otimes\mathcal V_{\mathbf t}$, it suffices to evaluate these maps on the highest vectors $e_{\mathbf s}$ and $e_{\mathbf t}$. By Theorem~\ref{thm:berezin-eigenvalues} and the preceding proposition,
\[
\sigma_m(T_m(e_{\mathbf s}))=\frac{A_{\mathbf s}(m)}{D_{\mathbf s}(m)}e_{\mathbf s}=\sigma_m\bigl(D_{\mathbf s}(m)^{-1}C_{\mathbf s}^{(m)}\bigr).
\]
Injectivity of $\sigma_m$ gives
\[
T_m(e_{\mathbf s})=D_{\mathbf s}(m)^{-1}C_{\mathbf s}^{(m)}.
\]
The identity $C_{\mathbf s}^{(m)}C_{\mathbf t}^{(m)}=C_{\mathbf s+\mathbf t}^{(m)}$ now proves the assertion.
\end{proof}
For $m\geq s_1$, define
\[
\mathcal Q_{m,\mathbf s}\coloneq D_{\mathbf s}(m)\,T_m|_{\mathcal V_{\mathbf s}}\colon\mathcal V_{\mathbf s}\longrightarrow\operatorname{End}(\mathcal H_m).
\]
Then
\begin{equation}\label{eq:pre-fixed-symbol}
\sigma_m(\mathcal Q_{m,\mathbf s}(f))=A_{\mathbf s}(m)f,\qquad \mathcal Q_{m,\mathbf s}(e_{\mathbf s})=C_{\mathbf s}^{(m)}.
\end{equation}
In particular, $\mathcal Q_{m,\mathbf s}$ is injective. Since $\sigma_m$ embeds $\operatorname{End}(\mathcal H_m)$ equivariantly into the multiplicity-free space of functions on $M$, the image of $\mathcal Q_{m,\mathbf s}$ is the entire isotypic component of type $\mathbf s$. Denote the projection onto this component by $\operatorname{pr}^{\operatorname{End}}_{\mathbf s}$. The preceding theorem gives, for $m\geq s_1+t_1$,
\begin{equation}\label{eq:pre-fixed-product}
\operatorname{pr}^{\operatorname{End}}_{\mathbf s+\mathbf t}\bigl(\mathcal Q_{m,\mathbf s}(f)\mathcal Q_{m,\mathbf t}(g)\bigr)=\mathcal Q_{m,\mathbf s+\mathbf t}(f\odot g).
\end{equation}
Iteration yields the following consequence.
\begin{corollary}\label{thm:interated-exact-cartan}
Let $f_i\in\mathcal V_{\mathbf s^{(i)}}$, $1\leq i\leq n$. For every integer $m\geq\sum_{i=1}^n s^{(i)}_1$,
\[
\operatorname{pr}_{\mathbf s^{(1)}+\cdots+\mathbf s^{(n)}}\sigma_m\bigl(T_m(f_1)\cdots T_m(f_n)\bigr)=\frac{A_{\mathbf s^{(1)}+\cdots+\mathbf s^{(n)}}(m)}{\prod_{i=1}^n D_{\mathbf s^{(i)}}(m)}\,(f_1\odot\cdots\odot f_n).
\]
\end{corollary}
Finally, we state the nonvanishing property of Cartan powers that will be used in Section~\ref{sec:nonlinear-obstructions}.
\begin{lemma}\label{lem:pre-cartan-nonvanishing}
If $0\neq x\in\mathcal V_{\mathbf s}$, then $x^{\odot d}\neq 0$ in $\mathcal V_{d\mathbf s}$ for every integer $d\geq1$.
\end{lemma}
\begin{proof}
Equip $\mathcal V_{\mathbf s}$ with a $G^*$-invariant Hermitian inner product and choose a nonzero highest vector $v$. The tensor $v^{\otimes d}$ lies in the irreducible component of $\mathcal V_{\mathbf s}^{\otimes d}$ with highest weight $d\Lambda_{\mathbf s}$, which is isomorphic to $\mathcal V_{d\mathbf s}$. If $x^{\odot d}=0$, then for every $u\in G^*$,
\[
0=\langle (u\cdot x)^{\otimes d},v^{\otimes d}\rangle=\langle u\cdot x,v\rangle^d.
\]
Since the $G^*$-orbit of $v$ spans the irreducible representation $\mathcal V_{\mathbf s}$, we must have $x=0$, a contradiction.
\end{proof}

\section{Variational formulae}\label{sec:variation}
We derive the first and second variational formulae for the canonically twisted Bergman kernel at the invariant metric $h_0$ in this section. The first variation is expressed in terms of the Berezin transform and the Laplacian and will be used to identify the infinitesimal deformations. We then relate these variations to the Bergman coefficients and show how poles of the resulting rational functions determine whether their first variations vanish.

\subsection{The varying Gram matrix}\label{subsec:varyingGram}
Let $(M,\omega_0;L,h_0)$ be the standard invariant polarization on an irreducible Hermitian symmetric space of compact type, as we described before. Set $h_t=h_0e^{-t\varphi}$ and $\omega_t=\omega_0+t\ii\pa\pp \varphi>0$, then
\begin{equation}\label{eq:volume-first-variation}
\left.\frac{d}{dt}\right |_{t=0}\log \det g_t=-\Delta \varphi,
\end{equation}
while the second differentiation is
\[
\left.\frac{d^2}{dt^2}\right |_{t=0}\log \det g_t=-|\ii\pa\pp \varphi|^2.
\]
For a fixed integer $k\geq p$, put
\[
m=k-p,\qquad L^k\otimes K_M\simeq L^m
\]
Locally, the metric on $K_M$ induced by $\omega_0$ is given by $(\det g_0)^{-1}$ and is invariant. After rescaling the isomorphism by a nonzero constant, we may therefore assume that it identifies the induced metric on $L^k\otimes K_M$ with $h_0^m$. Let $s_1,\ldots,s_{d_m}$ be an orthonormal basis of $\mathcal H_m$ with respect to the reference $L^2$ inner product, define the varying Gram matrix $G_t$ by
\[
(G_t)_{ij}=\langle s_i,s_j\rangle_t,\qquad 1\leq i,j\leq d_m,
\]
where $\langle\cdot,\cdot\rangle_t$ is the $L^2$ inner product induced by $h_t^k$ and $\omega_t$.

\begin{proposition}\label{prop:varying-identity}
Let $T_m$ and $\sigma_m$ denote the Toeplitz map and the covariant symbol map associated with $\mathcal H_m=H^0(M,L^m)$, then
\begin{equation}\label{eq:varyingGram}
G_t=T_m(e^{-kt\varphi})
\end{equation}
and
\begin{equation}\label{eq:varyingbergman}
\frac{\mathcal B^{K_M}_k(h_t;x)}{d_m}=q_t(x) \sigma_m(G_t^{-1})(x),\qquad q_t=e^{-kt\varphi}\frac{\det g_0}{\det g_t}.
\end{equation}
\end{proposition}
\begin{proof}
Relative to the initial metric, the pointwise metric on $L^k\otimes K_M$ is $e^{-kt\varphi}\frac{\det g_0}{\det g_t}$. The new probability measure $d\mu_t=\frac{\det g_t}{\det g_0}d\mu_0$ cancels the determinant ratio in the $L^2$ product. Therefore the matrix entry of $G_t$ is the matrix entry of the Toeplitz operator with symbol $e^{-kt\varphi}$, proving \eqref{eq:varyingGram}.

The sections $\{G_t^{-1/2}s_j\}_{j=1}^{d_m}$ form an orthonormal basis for the varying inner product. By the reproducing kernel identity,
\[
\begin{aligned}
\mathcal B^{K_M}_k(h_t;x)
&=q_t(x)\sum_{j=1}^{d_m}\|(G_t^{-1/2}s_j)(x)\|_{h_0^m}^2\\
&=q_t(x)\operatorname{tr}_{L_x^m}\bigl((G_t^{-1}K_m(\cdot,x))(x)\bigr)\\
&=q_t(x)d_m\,\sigma_m(G_t^{-1})(x),
\end{aligned}
\]
where we used $K_m(x,x)=d_m\operatorname{Id}_{L_x^m}$. This proves \eqref{eq:varyingbergman}.
\end{proof}

\subsection{The first and second variations}\label{subsec:first-second-variations}

\begin{theorem}\label{thm:first-variation}
For every smooth real function $\varphi$,
\begin{equation}\label{eq:first-variation}
\frac{1}{d_m}D\mathcal B^{K_M}_k|_{h_0}(\varphi)=\left( kB_m-k+\Delta\right)\varphi,\qquad m=k-p.
\end{equation}
\end{theorem}
\begin{proof}
From \eqref{eq:varyingGram}, 
\begin{equation}\label{eq:first-variation--Gram}
G'_0=-kT_m(\varphi),\qquad (G^{-1})'_0=kT_m(\varphi).
\end{equation}
Equations \eqref{eq:volume-first-variation} and \eqref{eq:varyingbergman} give
\begin{equation}\label{eq:first-variation--q}
q'_0=-k\varphi+\Delta \varphi.
\end{equation}
Differentiating \eqref{eq:varyingbergman} at $t=0$ and using $\sigma_m(I)=1$ and $\sigma_m(T_m(\varphi))=B_m(\varphi)$ yields \eqref{eq:first-variation}.
\end{proof}

Although the second variation is not used in the subsequent proof, we compute it for completeness. Define
\[
C_m(f,g)=\frac12\sigma_m\bigl(T_m(f)T_m(g)+T_m(g)T_m(f)\bigr).
\]

\begin{theorem}
\label{thm:second-variation}
For every smooth real function $\varphi$,
\begin{equation}\label{eq:second-variation}
\begin{aligned}
\frac{1}{d_m}D^2\mathcal B_k^{K_M}|_{h_0}(\varphi,\varphi)={}&(-k\varphi+\Delta\varphi)^2+|\ii\pa\pp\varphi|^2\\
&+2k(-k\varphi+\Delta\varphi)B_m(\varphi)\\
&+2k^2C_m(\varphi,\varphi)-k^2B_m(\varphi^2).
\end{aligned}
\end{equation}
\end{theorem}

\begin{proof}
Differentiating \eqref{eq:varyingGram} twice gives
\[
G''_0=k^2T_m(\varphi^2),\qquad (G_t^{-1})''_0=2k^2T_m(\varphi)^2-k^2T_m(\varphi^2).
\]
Moreover,
\[
q''_0=(-k\varphi+\Delta\varphi)^2+|\ii\pa\pp\varphi|^2.
\]
Differentiating \eqref{eq:varyingbergman} twice, using \eqref{eq:first-variation--Gram} and \eqref{eq:first-variation--q}, together with
$\sigma_m(T_m(\varphi^2))=B_m(\varphi^2)$ and $\sigma_m(T_m(\varphi)^2)=C_m(\varphi,\varphi)$, yields \eqref{eq:second-variation}.
\end{proof}

\begin{remark}
As shown in Proposition~\ref{prop:nonzero-density-pole}, the Taylor coefficient arising from our second variation suffices to establish the nonlinear obstruction when $a\geq4$ is even, whereas higher-order Taylor coefficients may be needed when $a=2$.
\end{remark}

\subsection{Variations of Bergman coefficients}
\label{subsec:differentiating-expansion}

To pass from the variational formulae above to the variations of Bergman coefficients, we need to differentiate the asymptotic expansion term by term. For the ordinary Bergman kernel, uniformity of the expansion and its compatibility with the first variation were established by Lu--Tian \cite[Proposition~4.1 and Theorems~4.1--4.2]{LuTian2004}. We use the following version for the canonically twisted kernel, including higher parameter derivatives.

\begin{lemma}
\label{lem:parameter-bergman}
Let $h_t$, $|t|\leq\delta$, be a smooth family of positive Hermitian metrics on $L$. For every $l,s,N\in\mathbb Z_{\geq 0}$, there exist constants $C_{l,s,N}>0$ and $k_0$ such that
\[
\sup_{|t|\leq\delta}\left\|\pa_t^l\left(\mathcal B_k^{K_M}(h_t)-\sum_{j=0}^{N}a^{K_M}_j(h_t)k^{n-j}\right)\right\|_{C^s}\leq C_{l,s,N}k^{n-N-1},\qquad k\geq k_0.
\]
The norms are defined using a fixed background metric on $M$.
\end{lemma}
\begin{proof}
Apply the parameter-uniform Bergman expansion of Dai--Liu--Ma \cite[Theorems~1.1 and~1.3]{DaiLiuMa2006} with $E=K_M$. The uniformity statement following their equation~(1.3) includes derivatives in the parameters. Smoothness and
compactness of the family give the required bounds.
\end{proof}

In particular, for every fixed $l$,
\begin{equation}\label{eq:differentiated-bergman}
\left.\pa_t^l\right|_{t=0}\mathcal B_k^{K_M}(h_t)\sim\sum_{j=0}^{\infty}\left.\pa_t^l\right|_{t=0}a^{K_M}_j(h_t)\,k^{n-j}.
\end{equation}
Certainly the same expansion holds after projection onto any fixed spherical summand. 

At the invariant metric, the first variation on each spherical module is described by an explicit rational function. Set
\[
\ell_{\mathbf s}(m)=(m+p)b_{\mathbf s}(m)-(m+p)+\lambda_{\mathbf s},
\]
By Theorem~\ref{thm:berezin-eigenvalues}, Theorem~\ref{thm:laplace-spectrum} and \eqref{eq:first-variation},
\begin{equation}\label{eq:spherical-scalar}
D\mathcal B_k^{K_M}|_{h_0}(\varphi)=d_m\ell_{\mathbf s}(m)\varphi,\qquad\text{for }\varphi\in\mathcal V_{\mathbf s}.
\end{equation}
Set
\[
R_{\mathbf s}(k)=d_{k-p}\ell_{\mathbf s}(k-p).
\]
By \eqref{eq:differentiated-bergman}, the coefficient of $k^{n-j}$ in the Laurent expansion of $R_{\mathbf s}(k)$ is the scalar by which $Da^{K_M}_j(0)$ acts on $\mathcal V_{\mathbf s}$.

\begin{lemma}\label{lem:rational-pole-principle}
Let $R(k)$ be a rational function with asymptotic expansion
\[
R(k)\sim\sum_{j=0}^{\infty}r_jk^{n-j}\qquad (k\to+\infty).
\]
Then $r_j=0$ for every $j\geq n+2$ if and only if
\begin{equation}\label{eq:rational-pole-principle}
R(k)=P(k)+\frac{c}{k},
\end{equation}
where $P$ is a polynomial of degree at most $n$ and $c\in\mathbb C$. In particular, such a function has no pole at any $k\neq 0$.
\end{lemma}

\begin{proof}
A rational function has a convergent Laurent expansion near infinity, whose coefficients agree with those of its asymptotic expansion. The assumed vanishing removes all powers $k^{-2},k^{-3},\ldots$. Hence \eqref{eq:rational-pole-principle} holds near infinity and therefore identically as an equality of rational functions. The converse is obvious.
\end{proof}
Consequently, if $R_{\mathbf s}(k)$ has a finite nonzero pole, then $Da^{K_M}_j(0)|_{\mathcal V_{\mathbf s}}\neq 0$ for some $j\geq n+2$.

\section{Classification of infinitesimal deformations}
\label{sec:infinitesimal-deformations}

In this section, we prove Theorem~\ref{thm:linearization} by examining each spherical summand separately. By \eqref{eq:spherical-scalar}, the first variation acts by a scalar on each summand, reducing the classification to an analysis of the corresponding rational functions.

For a partition $\mathbf s\in \mathcal P_r$, recall the numerator and denominator polynomials introduced in Theorem~\ref{thm:berezin-eigenvalues}
\[
A_{\mathbf s}(m)=(-1)^{|\mathbf s|}(-m)_{\mathbf s}=\prod_{j=1}^r\prod_{l=0}^{s_j-1}(m+c_j-l),
\]
and
\[
D_{\mathbf s}(m)=(m+p)_{\mathbf s}=\prod_{j=1}^r\prod_{l=0}^{s_j-1}(m+p-c_j+l),\qquad \text{where }c_j=\frac{a}{2}(j-1).
\]
Thus $b_{\mathbf s}(m)=A_{\mathbf s}(m)/D_{\mathbf s}(m)$.

\begin{lemma}\label{lem:no-cancellation}
The polynomials $A_{\mathbf s}(m)$ and $D_{\mathbf s}(m)$ have no common zero.
\end{lemma}
\begin{proof}
A zero of a denominator factor has the form 
\[
m=-p+c_j-l,\qquad 0\leq l<s_j.
\]
At this point any numerator factor in $(-m)_{\mathbf s}$ has value
\begin{equation}\label{eq:numerator-factor}
-m-c_i+u=p-c_j+l-c_i+u,\qquad 0\leq u<s_i.
\end{equation}
Since $c_i,c_j\leq a(r-1)/2$ and $p=2+a(r-1)+b$, every expression in \eqref{eq:numerator-factor} is therefore strictly positive.
\end{proof}

Now define
\[
\mathcal K_{\operatorname{lin}}\coloneq\bigcap_{j\geq n+2} \ker Da_j|_{h_0}\subset C^{\infty}(M,\mathbb{R})
\]
be the common kernel of the linearization of the system of Bergman logarithmic flatness equations.

\begin{proposition}\label{prop:divisibility-criterion}
The complex spherical summand $\mathcal V_{\mathbf s}$ is contained in the complexification of $\mathcal K_{\operatorname{lin}}$ if and only if
\begin{equation}\label{eq:divisibility-criterion}
D_{\mathbf s}(m) \mathrel{\big |} (m+p)^2d_m
\end{equation}
in $\mathbb{C}[m]$.
\end{proposition}
\begin{proof}
Put $k=m+p$, by Theorem~\ref{thm:first-variation}, the exact scalar on $\mathcal V_{\mathbf s}$ is
\[
R_{\mathbf s}(k)=d_{m}\ell_{\mathbf s}(m)=d_m\left(k\frac{A_{\mathbf s}(m)}{D_{\mathbf s}(m) } -k+\lambda_{\mathbf s}\right),
\]
thus $\mathcal V_{\mathbf s}$ lies in the complexification of $\mathcal K_{\operatorname{lin}}$ exactly when the rational function $R_{\mathbf s}(k)$ has the form $P(k)+\frac{c}{k}$. Equivalently, $kd_m\ell_{\mathbf s}(m)$ is a polynomial. 

In fact the dimension
\begin{equation}\label{eq:dimension-polynomial}
d_m=C_M\prod_{j=1}^{r}(m+1+c_j)_{1+b+a(r-j)}
\end{equation}
is a Hilbert polynomial, where $C_M$ is a constant such that $d_0=1$. The only possibly nonpolynomial term in $kd_m\ell_{\mathbf s}(m)$ is
\[
\frac{k^2d_mA_{\mathbf s}(m)}{D_{\mathbf s}(m)}.
\]
Then \eqref{eq:divisibility-criterion} follows from Lemma~\ref{lem:no-cancellation}.
\end{proof}

\begin{theorem}\label{thm:linearization-repeat}
If $a$ is odd, then
\begin{equation}\label{eq:odd-linearization}
\mathcal K_{\operatorname{lin}}=\mathbb{R}\mathbf{1}\oplus \mathcal{V}_{(1)}^{\mathbb{R}}.
\end{equation}
If $a$ is even, then
\begin{equation}\label{eq:even-linearization}
\mathcal K_{\operatorname{lin}}=\mathbb{R}\mathbf{1}\oplus \mathcal{V}_{(1)}^{\mathbb{R}}\oplus \bigoplus_{q=2}^r \mathcal{V}_{(1^q)}^{\mathbb{R}}.
\end{equation}
The first two summands are precisely the infinitesimal constant rescalings and the normalized real holomorphy potentials.
\end{theorem}
\begin{proof}
Suppose first $s_1\geq 2$. The factor $(m+p)_{s_1}$ in $D_{\mathbf s}(m)$ contains $m+p+1$, whose zero is $m=-p-1$. This zero is not supplied by $(m+p)^2$. By \eqref{eq:dimension-polynomial} the roots of $d_m$ are
\begin{equation}\label{eq:root-Hilbert}
m=-(1+c_j+u),\qquad 0\leq u\leq b+a(r-j).
\end{equation}
For every such root,
\[
1+c_i+u\leq 1+c_i+b+a(r-i)=p-1-c_i\leq p-1.
\]
Hence $-p-1$ is not a root of $d_m$. Criterion \eqref{eq:divisibility-criterion} fails. It follows that every possible partition has $s_1\leq 1$, and because the partition is decreasing, it is a column
\[
(1^q)=(\underbrace{1,\ldots,1}_{q},0,\ldots,0),\qquad 0\leq q\leq r.
\]
For such a column,
\begin{equation}\label{eq:factor-denominator}
D_{(1^q)}(m)=\prod_{j=1}^{q}(m+p-c_j).
\end{equation}
The $j=1$ factor is $m+p$ and is supplied by $(m+p)^2$.

Assume $a$ is even, then every $c_j$ is an integer. For $j\geq 2$, the root $m=-p+c_j$ occurs among the $j=1$ factors of $d_m$: in \eqref{eq:root-Hilbert} take $u=p-c_j-1$, this is an integer and $u\leq b+a(r-1)$. Thus all factors in \eqref{eq:factor-denominator} divide $(m+p)^2d_m$, this proves \eqref{eq:even-linearization}.

Assume $a$ is odd and $q\geq 2$. The $j=2$ factor of \eqref{eq:factor-denominator} vanishes at $m=-p+\frac{a}{2}$. To solve the equation $p-\frac{a}{2}=1+c_j+u$, firstly $j=2l$ must be even. This gives
\[
u-\left(b+a(r-j)\right)=1+a(l-1)>0
\]
This contradiction rules out every such column of height at least two. For $q=0$, $\ell_{(0)}(k-p)=k-k+0=0$; For $q=1$, $\ell_{(1)}(k-p)=k\frac{k-p}{k}-k+p=0$, which yields \eqref{eq:odd-linearization}. The last assertion follows from Proposition~\ref{prop:constants-holomorphy}.
\end{proof}

\section{Elliptic reduction modulo symmetries}
\label{sec:elliptic-reduction}

The main goal of this section is to prove Proposition~\ref{prop:finite-dimensional-graph}, which gives a finite-dimensional real-analytic parametrization of the solutions near the invariant metric. We first use the ellipticity of the linearized Bergman coefficient equations to construct a finite-order elliptic equation whose linearization has the same kernel as the original infinite system. After removing the directions corresponding to constant rescalings and holomorphic automorphisms, we apply the implicit-function theorem to show that every sufficiently small solution of \eqref{eq: logflatequation}, after normalization, lies on a finite-dimensional real-analytic graph parametrized by the remaining infinitesimal deformation directions. This reduction allows us to study the nonlinear obstructions in the subsequent sections.

\subsection{The coefficient maps and their linearizations}
\label{subsec:coefficient-maps}

For $\varphi\in C^{\infty}(M)$ define
\begin{equation}\label{eq:coefficient-equation}
\mathcal A_j(\varphi)=a^{K_M}_j(h_0e^{-\varphi}),\qquad\mathcal P_j=D\mathcal A_j(0),\qquad j\geq1.
\end{equation}
Since the reference metric $\omega_0$ is Bergman logarithmically flat,
\[
\mathcal A_j(0)=0,\qquad j\geq n+2.
\]

Lu \cite[Theorem~1.1]{Lu2000} expressed the ordinary Bergman coefficients as universal polynomials in the curvature and its covariant derivatives. Let $a_j$ denote the ordinary Bergman coefficients, normalized by $a_0=1$. Lu--Tian~\cite[Theorem~3.1]{LuTian2004} proved that
\[
a_j=\frac{(-1)^{j-1}j}{(j+1)!}\Delta^{j-1}\rho +\text{terms of lower differential order in the K\"ahler potential},
\]
where $\rho$ is the scalar curvature. Consequently, prescribing $a_j$ gives an elliptic equation of order $2j+2$ in the K\"ahler potential.

The following lemma gives the corresponding properties for the canonically twisted Bergman coefficients, together with their real-analytic dependence on the K\"ahler potential in H\"older spaces.

\begin{lemma}\label{lem:coefficient-analyticity}
For every $j\geq1$, every nonnegative integer $s$, and every $\alpha\in(0,1)$, set
\[
\mathcal U_{j,s}=\left\{\varphi\in C^{2j+2+s,\alpha}(M,\mathbb R)\colon\omega_\varphi>0\right\}.
\]
The coefficient map $\mathcal A_j$ extends to a real-analytic map
\[
\mathcal A_j\colon\mathcal U_{j,s}\longrightarrow C^{s,\alpha}(M,\mathbb R).
\]
Its derivative at zero is a bounded operator
\[
\mathcal P_j\colon C^{2j+2+s,\alpha}(M,\mathbb R)\longrightarrow C^{s,\alpha}(M,\mathbb R),
\]
which is a $G^*$-equivariant, formally self-adjoint elliptic differential operator of order $2j+2$. More precisely,
\begin{equation}\label{eq:linearized-coefficient-principal-part}
\mathcal P_j=\frac{(-1)^{j+1}C_M}{(j+1)!}\Delta^{j+1}+\mathcal E_j,\qquad \operatorname{ord}\mathcal E_j\leq 2j+1,
\end{equation}
where $C_M>0$ is the leading coefficient of the dimension polynomial $d_m$.
\end{lemma}

\begin{proof}
The Bergman kernel expansion with a Hermitian holomorphic vector bundle $E$ twisted expresses each coefficient as a universal polynomial in the curvatures and their covariant derivatives; see \cite[Theorems~1.1 and~1.3 and \S4]{DaiLiuMa2006}. Apply this by taking $E=K_M$ equipped with the metric induced by $\omega_\varphi$. The coefficient of index $j$ has weight $j$, where curvature has weight one and each covariant derivative has weight one half. Thus at most $2j-2$ covariant derivatives fall on any curvature factor. Both the curvature of $\omega_\varphi$ and that of $K_M$ involve at most four derivatives of $\varphi$. Consequently, $\mathcal A_j(\varphi)$ depends on derivatives of $\varphi$ of order at most $2j+2$.

In local coordinates, this dependence is obtained by differentiation, multiplication, and inversion of the positive matrix $g_\varphi$. These operations are real analytic on the corresponding Hölder spaces. This proves the asserted analyticity and boundedness of the derivative.

To identify the principal part, write the Berezin expansion at $h_0$ as
\[
B_m\sim\sum_{q=0}^{\infty}m^{-q}Q_q, \qquad Q_0=I.
\]
The principal part of $Q_j$ is $(-1)^j\Delta^j/j!$, as follows from \cite[equations~(2.23) and~(3.9)]{Englis2000}.

By \eqref{eq:first-variation} and the differentiated
Bergman expansion,
\[
d_m\bigl(kB_m-k+\Delta\bigr)\sim\sum_{j\geq 0}k^{n-j}\mathcal P_j,\qquad m=k-p.
\]
Since
\[
d_{k-p}=C_Mk^n+O(k^{n-1}),
\]
only $Q_q$ with $q\le j+1$ contribute to the coefficient of $k^{n-j}$. Among these, the only contribution of differential order $2j+2$ comes from the principal part of $Q_{j+1}$. This gives \eqref{eq:linearized-coefficient-principal-part}, so $\mathcal P_j$ is elliptic.

The coefficient construction and the invariance of $h_0$ imply $G^*$-equivariance of $\mathcal P_j$. Finally, $B_m$ is self-adjoint in $L^2(M,d\mu_0)$, hence $d_m\bigl(kB_m-k+\Delta\bigr)$ is self-adjoint for every sufficiently large integer $k$. Pairing its asymptotic expansion with arbitrary smooth functions and comparing coefficients gives
\[
\langle\mathcal P_jf,g\rangle_{L^2(d\mu_0)}=\langle f,\mathcal P_jg\rangle_{L^2(d\mu_0)}.
\]
Thus $\mathcal P_j$ is formally self-adjoint.
\end{proof}

\subsection{Reduction to a finite elliptic system}
\label{subsec:elliptic-reduction}
The first operator $\mathcal P_{n+2}$ is elliptic by Lemma~\ref{lem:coefficient-analyticity}, so its kernel is finite dimensional. We consider the following descending chain inside the kernel
\[
K_N=\bigcap_{j=n+2}^N \ker \mathcal P_j,\qquad N=n+2,n+3,\ldots .
\]
Every strict inclusion lowers the dimension, hence the chain stabilizes eventually. Choose $N_0\geq n+2$ such that
\begin{equation}\label{eq:actual-linearized-kernel}
K_{N_0}=\bigcap_{j=n+2}^{N_0} \ker \mathcal P_j=\mathcal K_{\operatorname{lin}}.
\end{equation}
Set
\[
d=4N_0+4,\qquad X=C^{d,\alpha}(M,\mathbb R),\qquad Y=C^{0,\alpha}(M,\mathbb R),
\]
and define the single finite-order map
\[
\mathcal F_{N_0}(\varphi)=\sum\limits_{j=n+2}^{N_0}\mathcal P_j\left(\mathcal A_j(\varphi)\right),\qquad \mathcal F\coloneq \mathcal F_{N_0}.
\]
Here each $\mathcal P_j$ is evaluated at the reference metric $h_0$. For every nonnegative integer $s$, the same expression defines
a real-analytic map
\[
\mathcal F\colon\mathcal U_{2N_0+1}=\{\varphi\in C^{d+s,\alpha}(M,\mathbb R)\colon\omega_\varphi>0\}\longrightarrow C^{s,\alpha}(M,\mathbb R).
\]
Indeed, for $n+2\leq j\leq N_0$, each composition factors as
\[
C^{d+s,\alpha}\hookrightarrow C^{4j+4+s,\alpha}\xrightarrow{\mathcal A_j} C^{2j+2+s,\alpha}\xrightarrow{\mathcal P_j} C^{s,\alpha},
\]
where $\mathcal A_j$ is restricted to the open set of positive potentials.

\begin{proposition}\label{prop:finite-elliptic-reduction}
The linearization of $\mathcal F$ at zero is
\[
\mathcal L\coloneq D\mathcal F(0)=\sum_{j=n+2}^{N_0}\mathcal P_j^2.
\]
It is a formally self-adjoint elliptic differential operator of order $d$, and
\[
\ker\mathcal L=\mathcal K_{\operatorname{lin}}.
\]
Moreover, every positive potential satisfying $\mathcal A_j(\varphi)=0$ for $n+2\leq j\leq N_0$ is a zero of $\mathcal F$.
\end{proposition}
\begin{proof}
Since $D\mathcal A_j(0)=\mathcal P_j$, differentiation gives the formula for $\mathcal L$. Its only summand of order $d=4N_0+4$ is
$\mathcal P_{N_0}^2$, whose principal symbol is the square of the nonvanishing principal symbol of $\mathcal P_{N_0}$. Thus $\mathcal L$ is elliptic. Formal self-adjointness follows from that of each $\mathcal P_j$.

If $u\in X$ satisfies $\mathcal Lu=0$, elliptic regularity then implies that $u$ is smooth. Integration by parts then gives
\[
0=\langle\mathcal Lu,u\rangle_{L^2(d\mu_0)}=\sum_{j=n+2}^{N_0}\|\mathcal P_ju\|_{L^2(d\mu_0)}^2.
\]
Hence $u\in K_{N_0}=\mathcal K_{\operatorname{lin}}$. The reverse inclusion is immediate. The last assertion follows directly from the definition of $\mathcal F$.
\end{proof}

By elliptic regularity and the Fredholm alternative, $\mathcal L\colon X\longrightarrow Y$ is Fredholm of index zero, and
\[
\operatorname{Ran}\mathcal L=\left\{f\in Y\colon\langle f,u\rangle_{L^2(d\mu_0)}=0\text{ for every }u\in\mathcal K_{\operatorname{lin}}\right\}.
\]
The global Schauder estimate gives
\[
\|u\|_{C^{d,\alpha}}\leq C\left(\|\mathcal Lu\|_{C^{0,\alpha}}+\|u\|_{C^0}\right),\qquad u\in X.
\]

\subsection{Removing the symmetry directions}
\label{subsec:symmetry-directions}
Set $\mathcal S=\mathbb R\mathbf 1\oplus\mathcal V_{(1)}^{\mathbb R}$. All orthogonal complements and projections below are taken with respect to $L^2(M,d\mu_0)$. For a finite-dimensional space $W$ of smooth functions, let $\Pi_W$ denote the orthogonal projection onto $W$. 

Choose a real vector-space complement $\mathfrak q$ to $\mathfrak g^*$ in the Lie algebra of $\operatorname{Aut}^0(M,L)$. For $A\in\mathfrak q$ sufficiently small, write $F_A=\exp A$, using the same notation for its holomorphic lift to $L$. Define $\Gamma(A,c)$ by
\[
e^{-c}F_A^*h_0=h_0e^{-\Gamma(A,c)},\qquad c\in\mathbb R.
\]
We normalize the lift so that
\[
\int_M\Gamma(A,0)\,d\mu_0=0.
\]
Then $\Gamma(A,c)=\Gamma(A,0)+c$. For $A\in\mathfrak q$, the function $u_A\coloneq\left.\frac{d}{dt}\right|_{t=0}\Gamma(tA,0)$ is the normalized potential of the infinitesimal deformation induced by $F_{tA}^*h_0$. The map $A\mapsto u_A$ is a linear isomorphism from $\mathfrak q$ onto $\mathcal V_{(1)}^{\mathbb R}$. Consequently,
\begin{equation}\label{eq:symmetry-orbit-differential}
D\Gamma(0,0)\colon\mathfrak q\oplus\mathbb R\longrightarrow\mathcal S
\end{equation}
is an isomorphism. The Bergman coefficients $a^{K_M}_j(e^{-c}F^*h)=F^*a^{K_M}_j(h)$. Since the reference metric satisfies $a^{K_M}_j(h_0)=0$ for $j\geq n+2$, every potential $\Gamma(A,c)$ satisfies these equations as well. We next normalize the potential by requiring it to be $L^2$-orthogonal to the infinitesimal symmetry directions. 

\begin{lemma}\label{lem:local-moment-slice}
Let $s\geq2$ be an integer. For every sufficiently small $\varphi\in C^{s,\alpha}(M,\mathbb R)$, there is a unique sufficiently small pair
$(A,c)\in\mathfrak q\oplus\mathbb R$ such that the potential $\psi$ defined by
\[
e^{-c}F_A^*(h_0e^{-\varphi})=h_0e^{-\psi}
\]
satisfies
\[
\Pi_{\mathcal S}\psi=0.
\]
Moreover, $\psi\rightarrow 0$ in $C^{s,\alpha}$ as $\varphi\rightarrow 0$ in $C^{s,\alpha}$. For smooth potentials, this transformation preserves positivity of the K\"ahler form, Bergman logarithmic flatness, and each equation $a^{K_M}_j(h_\varphi)=0$.
\end{lemma}

\begin{proof}
Define
\[
\Psi(A,c,\varphi)=\Gamma(A,c)+\varphi\circ F_A,\qquad \mathcal M(A,c,\varphi)=\Pi_{\mathcal S}\Psi(A,c,\varphi).
\]
We apply the implicit-function theorem to $\mathcal M$. Since $s\geq 2$, the smooth dependence of $F_A$ on $A$ and the finite dimensionality of $\mathcal S$ imply that $\mathcal M$ is a $C^1$ map from a neighborhood of $(0,0,0)$ in $\mathfrak q\times\mathbb R\times C^{s,\alpha}(M,\mathbb R)$ to $\mathcal S$. At the origin,
\[
D_{(A,c)}\mathcal M(0,0,0)=D\Gamma(0,0),
\]
which is an isomorphism by \eqref{eq:symmetry-orbit-differential}. The implicit-function theorem therefore gives unique small $A(\varphi)$ and $c(\varphi)$ satisfying $\mathcal M=0$.

Set
\[
\psi=\Gamma(A(\varphi),c(\varphi))+\varphi\circ F_{A(\varphi)}.
\]
Since $A(\varphi),c(\varphi)\rightarrow 0$ as $\varphi\rightarrow0$ we obtain $\psi\rightarrow 0$ in $C^{s,\alpha}$. Finally, $\omega_\psi=F_{A(\varphi)}^*\omega_\varphi$, so positivity is preserved. The remaining assertions are obvious now.
\end{proof}

Define $E=\mathcal K_{\operatorname{lin}}\cap\mathcal S^\perp$. Then $\mathcal K_{\operatorname{lin}}=\mathcal S\oplus E$, where $E$ is a finite-dimensional $G^*$-module consisting of smooth functions. To remove the infinitesimal symmetry directions, we restrict the potentials to
\[
X_{\mathrm{sl}}\coloneq X\cap\mathcal S^\perp.
\]
Since $\mathcal K_{\operatorname{lin}}=\mathcal S\oplus E$, this space decomposes as
\[
X_{\mathrm{sl}}=E\oplus X_1,\qquad \text{where }X_1\coloneq X\cap\mathcal K_{\operatorname{lin}}^\perp.
\]
We also set
\[
Y_1\coloneq Y\cap\mathcal K_{\operatorname{lin}}^\perp.
\]
Then $\mathcal L:X_1\longrightarrow Y_1$ is a bounded isomorphism.

The particular construction of $\mathcal F$ gives the stronger identity
\begin{equation}\label{eq:finite-map-orthogonality}
\mathcal F(\varphi)
\perp\mathcal K_{\operatorname{lin}}
\end{equation}
for every $\varphi$ in its domain.
Indeed, if $u\in\mathcal K_{\operatorname{lin}}$, then $\mathcal P_j u=0$ for $n+2\leq j\leq N_0$. Formal self-adjointness therefore gives
\[
\begin{aligned}
\langle\mathcal F(\varphi),u\rangle_{L^2}&=\sum_{j=n+2}^{N_0}\langle\mathcal P_j\left(\mathcal A_j(\varphi)\right),u\rangle_{L^2}\\
&=\sum_{j=n+2}^{N_0}\langle\mathcal A_j(\varphi),\mathcal P_j u\rangle_{L^2}=0.
\end{aligned}
\]
In particular, $\mathcal F$ takes values in $Y_1$.

Consider the equation
\[
\mathcal F(v+u)=0,\qquad v\in E,\quad u\in X_1.
\]
Since $\mathcal F(0)=0$ and its derivative in $u$ at the origin is the isomorphism $\mathcal L\colon X_1\to Y_1$, the analytic implicit-function theorem gives a unique local real-analytic map
\[
w\colon U_E\longrightarrow X_1,\qquad w(0)=0,\qquad Dw(0)=0,
\]
where $U_E$ is a neighborhood of zero in $E$, such that
\[
\mathcal F(v+w(v))=0.
\]
Here $Dw(0)=0$ follows by differentiating this identity and using $\mathcal L|_E=0$.

\begin{proposition}\label{prop:finite-dimensional-graph}
Every sufficiently small $\varphi\in X_{\mathrm{sl}}$ satisfies $\mathcal F(\varphi)=0$ if and only if
\[
\varphi=v+w(v)
\]
for a unique sufficiently small $v\in E$. In particular, every sufficiently small solution $\varphi\in X_{\mathrm{sl}}$ of
\[
\mathcal A_j(\varphi)=0,\qquad n+2\leq j\leq N_0,
\]
lies on this graph.
\end{proposition}

\begin{proof}
Write $\varphi=v+u$ according to $X_{\mathrm{sl}}=E\oplus X_1$. Uniqueness in the implicit-function theorem gives $\mathcal F(v+u)=0$ if and only if $u=w(v)$. The last assertion follows from the definition of $\mathcal F$.
\end{proof}
It remains to determine which functions $\varphi=v+w(v)$ satisfy the original coefficient equations. Since $\mathcal F$ and the decomposition $X_{\mathrm{sl}}=E\oplus X_1$ are $G^*$-equivariant, uniqueness in the implicit-function theorem then implies that the map $w$ is also $G^*$-equivariant.

\subsection{Higher regularity and compatibility}
\label{subsec:regularity-compatibility}

Later we shall differentiate the Bergman kernel along families of potentials arising from the graph $\varphi=v+w(v)$. We therefore establish the regularity of the Taylor coefficients of $w$ needed for these calculations.

For every nonnegative integer $s$, set
\[
X_1^{(s)}=C^{d+s,\alpha}(M,\mathbb R)\cap\mathcal K_{\operatorname{lin}}^\perp,\qquad Y_1^{(s)}=C^{s,\alpha}(M,\mathbb R)\cap\mathcal K_{\operatorname{lin}}^\perp.
\]
Thus $X_1^{(0)}=X_1$ and $Y_1^{(0)}=Y_1$. The same differential expression defines a real-analytic map near zero,
\[
\mathcal F\colon E\oplus X_1^{(s)} \longrightarrow Y_1^{(s)}.
\]
The target follows from \eqref{eq:finite-map-orthogonality}. Moreover, elliptic regularity and the range characterization of $\mathcal L$ give a bounded isomorphism
\[
\mathcal L\colon X_1^{(s)} \longrightarrow Y_1^{(s)}.
\]

\begin{lemma}\label{lem:compatible-graph-germs}
For every nonnegative integer $s$, there is an open neighborhood $U_s\subset U_E$ of $0$ such that
\[
w|_{U_s}\colon U_s\longrightarrow X_1^{(s)}
\]
is real analytic. The neighborhoods may be chosen so that $U_{s+1}\subset U_s$. Consequently, every Taylor tensor $D^l w(0)$ has smooth image and is $G^*$-equivariant.
\end{lemma}

\begin{proof}
At regularity level $s$, the analytic implicit-function theorem gives a local real-analytic map
\[
\widetilde w_s\colon\widetilde U_s\longrightarrow X_1^{(s)}
\]
with $\widetilde w_s(0)=0$ and $\mathcal F(v+\widetilde w_s(v))=0$. Under the continuous inclusion $X_1^{(s)}\hookrightarrow X_1$, this map also solves the equation at the initial regularity level. Uniqueness therefore gives $\widetilde w_s=w$ on a sufficiently small neighborhood $U_s$ of $0$. Shrinking these neighborhoods successively, we may arrange that $U_{s+1}\subset U_s$.

For every $l\geq1$ and $v_1,\ldots,v_l\in E$, the value $D^l w(0)[v_1,\ldots,v_l]$ thus belongs to $C^{d+s,\alpha}(M,\mathbb R)$ for every $s$, and hence is smooth. Equivariance follows by differentiating the $G^*$-equivariance identity for $w$.
\end{proof}

In subsequent applications, we will work on $U_s$ for a fixed sufficiently large $s$.

\section{Local rigidity when $a$ is odd}
\label{sec:odd-a-rigidity}

We now prove Theorem~\ref{thm:localrigidity} when $a$ is odd. We restate this case below for the reader's convenience. The proof uses only the identity $\mathcal K_{\operatorname{lin}}=\mathbb R\mathbf1\oplus\mathcal V_{(1)}^{\mathbb R}$. It therefore applies whenever this kernel identity holds, including the rank-one case.
 
\begin{theorem}\label{thm:odd-a-rigidity}
Let $M$ be an irreducible Hermitian symmetric space of compact type of complex dimension $n$, equipped with its standard polarization $(L,h_0)$. Assume that the restricted-root parameter $a$ is odd. For every $\alpha\in(0,1)$, there exist integers $J\geq n+2$ and $s\geq 2$, and a constant $\varepsilon>0$, such that the following holds.

If $\varphi\in C^\infty(M,\mathbb R)$ satisfies
\[
h_\varphi=h_0e^{-\varphi},\qquad\omega_\varphi>0,\qquad\|\varphi\|_{C^{s,\alpha}(M,\omega_0)}<\varepsilon,
\]
and
\[
a^{K_M}_j(h_\varphi)\equiv0\qquad(n+2\leq j\leq J),
\]
then
\[
h_\varphi=CF^*h_0
\]
for some $C>0$ and $F\in\operatorname{Aut}^0(M,L)$. In particular, the conclusion holds if the corresponding unit circle bundle is Bergman logarithmically flat.
\end{theorem}

\begin{proof}
Fix $\alpha\in(0,1)$. By \eqref{eq:odd-linearization},
\[
\mathcal K_{\operatorname{lin}}=\mathbb R\mathbf1\oplus\mathcal V_{(1)}^{\mathbb R}.
\]
Thus $E=\{0\}$ and $X_{\mathrm{sl}}=X_1$. The graph constructed in Section~\ref{sec:elliptic-reduction} therefore consists only of the origin. Consequently, zero is the only sufficiently small solution of
\[
\mathcal F(\psi)=0,\qquad \psi\in X_{\mathrm{sl}}.
\]

Let $N_0$ be chosen as in \eqref{eq:actual-linearized-kernel}, and take
\[
J=N_0,\qquad s=d=4N_0+4.
\]
Suppose that $\varphi$ satisfies the hypotheses of the theorem. By Lemma~\ref{lem:local-moment-slice}, there exists a polarized automorphism $F_A$ near the identity and a real constant $c$ such that
\[
h_\psi=e^{-c}F_A^*h_\varphi,\qquad \psi\in X_{\mathrm{sl}}.
\]
By decreasing $\varepsilon$, we may ensure that $\psi$ lies in the neighborhood where the preceding uniqueness statement applies. Thus $\mathcal A_j(\psi)=F_A^*\mathcal A_j(\varphi)=0$ for all $n+2\leq j\leq N_0$.
Hence
\[
\mathcal F(\psi)=\sum_{j=n+2}^{N_0}\mathcal P_j\bigl(\mathcal A_j(\psi)\bigr)=0.
\]
It follows that $\psi=0$, and therefore
\[
h_\varphi=e^c(F_A^{-1})^*h_0.
\]
The final statement follows from the characterization of Bergman logarithmic flatness.
\end{proof}

\section{Nonlinear obstructions}
\label{sec:nonlinear-obstructions}

Throughout this section, assume that $a$ is even and $r\geq 2$. For $2\leq q\leq r$, set $\Lambda_q=\Lambda_{(1^q)}=\gamma_1+\cdots+\gamma_q$. The additional infinitesimal directions form $E=\bigoplus_{q=2}^r\mathcal V_{(1^q)}^{\mathbb R}$, with $E^{\mathbb C}=\bigoplus_{q=2}^r\mathcal V_{(1^q)}$.

We study the nonlinear coefficient equations on the finite-dimensional graph $\varphi=v+w(v)$, where $v\in E$. Our goal is to show that the origin is an isolated solution. To this end, we consider a nonconstant real-analytic curve through the origin and use highest-weight arguments to isolate the leading contribution of a single spherical component along this curve. The inverse Gram matrix recursion and the exact Cartan product identities then show that a suitable Taylor coefficient of the Bergman kernel has a rational dependence on $k$ with a nonzero pole at $k=a/2$. In the next section, we show that the vanishing of finitely many Bergman coefficients rules out such a pole, thereby establishing local rigidity.

\subsection{Highest-weight analysis along analytic curves}
\label{subsec:highest-weight-projections}

We write $\nu\preceq\lambda$ if $\lambda-\nu$ is a nonnegative integral combination of simple roots. Let $\operatorname{pr}_\lambda$ and $\operatorname{pr}^{\operatorname{End}}_\lambda$ denote the projections onto the isotypic components of highest weight $\lambda$ in the space of functions and in $\operatorname{End}(\mathcal H_m)$, respectively.

Suppose that $A$ and $B$ are finite sums of irreducible components with highest weights bounded by $i\lambda$ and $j\lambda$, respectively. The highest-weight theorem implies that $AB$ has highest weights bounded by $(i+j)\lambda$ and that
\begin{equation}\label{eq:highest-product-projection}
\operatorname{pr}_{(i+j)\lambda}(AB)=\operatorname{pr}_{(i+j)\lambda}\bigl((\operatorname{pr}_{i\lambda}A)(\operatorname{pr}_{j\lambda}B)\bigr).
\end{equation}
The same identity holds for operator composition, with $\operatorname{pr}_\lambda$ replaced by $\operatorname{pr}^{\operatorname{End}}_\lambda$.

Now we apply this observation to the graph map $w$. For $x\in\mathcal V_{(1^q)}$, consider the formal series
\[
\varphi_x(s)=sx+\sum_{d\geq 2}\frac{s^d}{d!}D^dw(0)[x,\ldots,x],
\]
where the Taylor tensors of $w$ are extended by complex multilinearity. For $m=k-p$, write
\[
G(s)=T_m(e^{-k\varphi_x(s)})=I+\sum_{d\geq 1}s^dG_d,\qquad G(s)^{-1}=I+\sum_{d\geq1}s^dH_d.
\]

\begin{lemma}\label{lem:highest-weight-support}
The coefficient of degree $d$ in each of
\[
\varphi_x(s),\quad e^{-k\varphi_x(s)},\quad q_{\varphi_x(s)},\quad G(s),\quad G(s)^{-1},\quad \sigma_m(G(s)^{-1})
\]
contains only irreducible components of highest weight $\nu\preceq d\Lambda_q$. Moreover,
\begin{equation}\label{eq:projected-inverse-recursion}
\operatorname{pr}^{\operatorname{End}}_{d\Lambda_q}H_d=-\operatorname{pr}^{\operatorname{End}}_{d\Lambda_q}G_d-\sum_{i=1}^{d-1}\operatorname{pr}^{\operatorname{End}}_{d\Lambda_q}\left((\operatorname{pr}^{\operatorname{End}}_{i\Lambda_q}G_i)(\operatorname{pr}^{\operatorname{End}}_{(d-i)\Lambda_q}H_{d-i})\right).
\end{equation}
\end{lemma}

\begin{proof}
Write
\[
\varphi_x(s)=\sum_{j\geq1}s^j\varphi_j(x),\qquad \varphi_1(x)=x,\qquad \varphi_j(x)=\frac1{j!}D^jw(0)[x,\ldots,x]\quad(j\geq 2).
\]
Equivariance of $w$ implies that each $\varphi_j$ is an equivariant homogeneous polynomial map of degree $j$, and hence induces an equivariant linear map from $\operatorname{Sym}^j\mathcal V_{(1^q)}$. Every irreducible summand of this symmetric power has highest weight at most $j\Lambda_q$, so the same bound holds for $\varphi_j(x)$.
For a formal power series $F(s)$, let $[s^d]F(s)$ denote the coefficient of $s^d$. For the exponential, we have
\[
[s^d]e^{-k\varphi_x(s)}=\sum_{l=1}^{d}\frac{(-k)^l}{l!}\sum_{\substack{j_1+\cdots+j_l=d\\ j_i\geq1}}\varphi_{j_1}(x)\cdots\varphi_{j_l}(x).
\]
Each term is homogeneous of degree $d$ in $x$, and pointwise multiplication is equivariant. Thus this coefficient is an equivariant homogeneous polynomial map of degree $d$.
For the determinant ratio, let $S(\varphi)$ be the endomorphism of $T^{1,0}M$ determined by
\[
\omega_0(S(\varphi)U,\overline V)=(\sqrt{-1}\pa\pp\varphi)(U,\overline V).
\]
Then
\[
\frac{\det g_0}{\det g_\varphi}=\det(I+S(\varphi))^{-1}.
\]
The map $S$ is linear and equivariant because $G^*$ acts by holomorphic isometries of $g_0$. Write the Taylor expansion at the zero endomorphism as
\[
\det(I+S)^{-1}=\sum_{l\geq0}P_l(S),
\]
where $P_l$ is a scalar-valued homogeneous polynomial of degree $l$ on $\operatorname{End}(T^{1,0}M)$. Substituting
\[
S(\varphi_x(s))=\sum_{j\geq1}s^jS(\varphi_j(x))
\]
shows that the coefficient of $s^d$ is an equivariant homogeneous polynomial of degree $d$ in $x$. The same is therefore true for the coefficient of $s^d$ in $q_{\varphi_x(s)}$.
Each of these degree-$d$ homogeneous polynomial maps determines a unique symmetric $d$-linear map, and hence an equivariant linear map from $\operatorname{Sym}^d\mathcal V_{(1^q)}$. Hence their images contain only irreducible summands of highest weight at most $d\Lambda_q$.
Since $G_d=T_m\bigl([s^d]e^{-k\varphi_x(s)}\bigr)$, equivariance of $T_m$ gives the same bound for $G_d$. Comparing coefficients in $G(s)G(s)^{-1}=I$ gives
\[
H_d=-G_d-\sum_{i=1}^{d-1}G_iH_{d-i}.
\]
Starting with $H_1=-G_1$, induction and the highest-weight bound for products give the assertion for every $H_d$. Applying \eqref{eq:highest-product-projection} to each product in this recursion yields \eqref{eq:projected-inverse-recursion}. Finally, equivariance of $\sigma_m$ gives the same bound for $\sigma_m(H_d)$.
\end{proof}

We evaluate the projected products in \eqref{eq:projected-inverse-recursion} using the generalized Cartan product normalization fixed in Subsection~\ref{subsec:pre-conical-cartan}. We next consider a nonconstant real-analytic curve in $E$ through the origin and show that the leading term of a suitable spherical component determines the corresponding highest-weight projections of the nonlinear expressions along the curve. This reduces the obstruction argument to a single spherical component.

Let $v(t)$ be a nonconstant real-analytic curve in $E$ with $v(0)=0$, and write
\[
v(t)=\sum_{j=2}^r v_j(t),\qquad \text{where }v_j(t)\in\mathcal V_{(1^j)}^{\mathbb R}.
\]
Set $I=\{j\colon v_j(t)\not\equiv 0\}$. For $j\in I$, write
\[
v_j(t)=\sum_{n\geq e_j}t^n x_{j,n},\qquad \text{where }x_j\coloneq x_{j,e_j}\neq 0 \text{ and } e_j\geq 1.
\]
Let $\Pi$ be the chosen set of simple roots and put
\[
\mathfrak h_{\mathbb R}=\operatorname{span}_{\mathbb R}\{h_\alpha\colon\alpha\in\Pi\}.
\]
Choose one $H\in\mathfrak h_{\mathbb R}$ such that $\alpha(H)>0$ $(\alpha\in\Pi)$ and the numbers
\[
\rho_j=\frac{\Lambda_j(H)}{e_j},\qquad j\in I,
\]
are pairwise distinct.
Such a choice is possible: the strongly orthogonal roots $\gamma_1,\ldots,\gamma_r$ are linearly independent, so the weights $\Lambda_j/e_j$ are distinct, and equality of two evaluations excludes only a proper hyperplane. Let $q\in I$ be the unique index for which $\rho_q=\max_{j\in I}\rho_j$. Since each $\Lambda_j$ is a sum of positive roots, $\rho_q>0$.

Consider a formal series
\[
\Psi(v)=\sum_{l\geq 0}\Psi_l(v,\ldots,v),
\]
where each $\Psi_l$ is a symmetric $G^*$-equivariant complex multilinear map
\[
\Psi_l\colon (E^{\mathbb C})^{l}\longrightarrow C^\infty(M,\mathbb C).
\]
For an analytic map $\Psi$, these are its Taylor coefficients: $\Psi_l=\frac1{l!}D^l\Psi(0)$. We use $[t^N]$ to denote the coefficient of $t^N$.

\begin{lemma}\label{lem:leading-highest-weight-coefficient}
For every $d\geq 1$,
\begin{equation}\label{eq:leading-highest-weight-coefficient}
[t^{de_q}]\operatorname{pr}_{d\Lambda_q}\Psi(v(t))=\operatorname{pr}_{d\Lambda_q}\Psi_d(x_q,\ldots,x_q)=[s^d]\operatorname{pr}_{d\Lambda_q}\Psi(sx_q).
\end{equation}
\end{lemma}

\begin{proof}
A term contributing to the coefficient of $t^{de_q}$
has the form
\[
\Psi_l(x_{q_1,n_1},\ldots,x_{q_l,n_l}),\qquad \text{where }q_i\in I,\quad n_i\geq e_{q_i}\text{ and } \sum_{i=1}^{l}n_i=de_q.
\]
If its projection onto the component of highest weight $d\Lambda_q$ is nonzero, equivariance and the highest-weight theorem imply that $d\Lambda_q\preceq \sum_{i=1}^{l}\Lambda_{q_i}$. Evaluating on $H$ gives
\[
\begin{aligned}
d\Lambda_q(H)&\leq \sum_{i=1}^{l}\Lambda_{q_i}(H)\\
&\leq \rho_q\sum_{i=1}^{l}e_{q_i}\\
&\leq \rho_q\sum_{i=1}^{l}n_i=d\Lambda_q(H).
\end{aligned}
\]
All inequalities are therefore equalities. Since $\rho_q$ is the unique maximum, the second equality forces $q_i=q$ for every $i$. Since $\rho_q>0$, the third equality gives $n_i=e_q$. Thus $l e_q=de_q$, which implies that $l=d$, and every $x_{q,e_q}=x_q$. This proves \eqref{eq:leading-highest-weight-coefficient}.
\end{proof}

Since $v(0)=0$, each term $\Psi_l(v(t),\ldots,v(t))$ vanishes to order at least $l$ at $t=0$. Thus $[t^N]\Psi(v(t))$ depends only on the Taylor coefficients of $v$ up to order $N$ and the maps $\Psi_l$ with $l\leq N$. The lemma therefore applies to the compatible Taylor
tensors of $w$, whose values are smooth, and to the formal expressions constructed from them.

\subsection{The pole obstruction}
\label{subsec:pole-obstruction}

Assume that $a$ is even and $r\geq 2$, and put
\[
A=\frac a2,\qquad c_i=A(i-1),\qquad m_0=A-p.
\]
We study the rational coefficients at $k=A$, where $m=k-p=m_0$. For a vector-valued rational function $R(k)$, let $\operatorname{pole}_{k=A}R$ denote its pole order at $A$, and $[R]_{-j}$ the coefficient of $(k-A)^{-j}$ in its Laurent expansion.

\begin{lemma}\label{lem:polynomial-zeros}
The polynomial $d_m$ has a simple zero at $m=m_0$. Moreover, $A_{\mathbf s}(m_0)\neq 0$ for every partition $\mathbf s$. For $2\leq q\leq r$ and $d\geq1$,
\begin{equation}\label{eq:denominator-zero-order}
\nu_d(q)\coloneq\operatorname{ord}_{k=A}D_{d(1^q)}(k-p)=\min\left\{q-1,\,1+\left\lfloor\frac{d-1}{A}\right\rfloor\right\},
\end{equation}
where $\operatorname{ord}$ denotes the multiplicity of zeros.
\end{lemma}

\begin{proof}
The Hilbert polynomial is
\[
d_m=C_M\prod_{j=1}^{r}(m+1+c_j)_{1+b+a(r-j)}
\]
for some $C_M>0$. A numerator factor vanishes at $m=m_0$ precisely when
\[
u=p-A-1-c_i,\qquad 0\leq u\leq b+2A(r-i).
\]
For $i=1$, this value lies in the indicated interval. For $i\geq  2$, it exceeds the upper endpoint by
\[
1+A(i-2)>0.
\]
Thus exactly one linear factor vanishes.

Every factor of $A_{\mathbf s}(m_0)$ has the form $m_0+c_i-u$, with $0\leq u<s_i$. Since
\[
-(m_0+c_i-u)=p-A-c_i+u\geq p-A-c_r=2+b+A(r-2)>0,
\]
none of these factors vanishes.

Finally,
\[
D_{d(1^q)}(k-p)=\prod_{i=1}^{q}\prod_{u=0}^{d-1}\bigl(k-A(i-1)+u\bigr).
\]
At $k=A$, a factor vanishes exactly when $u=A(i-2)$. The admissible indices are therefore
\[
2\leq i\leq \min\left\{q,\,2+\left\lfloor\frac{d-1}{A}\right\rfloor\right\},
\]
which proves \eqref{eq:denominator-zero-order}.
\end{proof}

Fix $2\leq q\leq r$ and $0\neq x\in\mathcal V_{(1^q)}$. Recall the formal series
\[
\varphi_x(s)=sx+w(sx),\qquad G(s)=T_m(e^{-k\varphi_x(s)})=I+\sum_{d\geq1}s^dG_d,\qquad G(s)^{-1}=I+\sum_{d\geq1}s^dH_d,
\]
where $m=k-p$. For integers $m\geq d$, the map $\mathcal Q_{m,d(1^q)}$ identifies $\mathcal V_{d(1^q)}$ with the entire corresponding isotypic component of $\operatorname{End}(\mathcal H_m)$. Thus there are unique $\gamma_d,\eta_d\in\mathcal V_{d(1^q)}$ such that
\[
\operatorname{pr}^{\operatorname{End}}_{d\Lambda_q}G_d=\mathcal Q_{m,d(1^q)}(\gamma_d),\qquad\operatorname{pr}^{\operatorname{End}}_{d\Lambda_q}H_d=\mathcal Q_{m,d(1^q)}(\eta_d).
\]
The vectors $\gamma_d$ and $\eta_d$ depend on $x$ and $k$.

\begin{lemma}\label{lem:projected-gram-coefficients}
These coefficients extend to rational functions of $k$ and satisfy
\begin{equation}\label{eq:gamma-explicit}
\gamma_d(k)=\frac{\operatorname{pr}_{d\Lambda_q}[s^d]e^{-k\varphi_x(s)}}{D_{d(1^q)}(k-p)},
\end{equation}
and
\begin{equation}\label{eq:eta-recursion}
\eta_d=-\gamma_d-\sum_{i=1}^{d-1}\gamma_i\odot\eta_{d-i}.
\end{equation}
In particular,
\[
\operatorname{pole}_{k=A}\gamma_d\leq\nu_d(q),\qquad \gamma_1(k)=-\frac{k}{D_{(1^q)}(k-p)}x,
\]
and $\gamma_1$ has a simple pole at $k=A$ with nonzero coefficient $[\gamma_1]_{-1}$.
\end{lemma}

\begin{proof}
By equivariance of $T_m$ and its relation to $\mathcal Q_{m,\mathbf s}$,
\[
\operatorname{pr}^{\operatorname{End}}_{d\Lambda_q}G_d=T_m\bigl(\operatorname{pr}_{d\Lambda_q}[s^d]e^{-k\varphi_x(s)}\bigr)
=\mathcal Q_{m,d(1^q)}\left(\frac{\operatorname{pr}_{d\Lambda_q}[s^d]e^{-k\varphi_x(s)}}{D_{d(1^q)}(k-p)}\right).
\]
Injectivity gives \eqref{eq:gamma-explicit}. Its numerator is polynomial in $k$. Thus Lemma~\ref{lem:polynomial-zeros} gives the pole bound.

Projecting the identity
\[
H_d=-G_d-\sum_{i=1}^{d-1}G_iH_{d-i}
\]
and using \eqref{eq:projected-inverse-recursion} and \eqref{eq:pre-fixed-product} yields \eqref{eq:eta-recursion}. This also proves rationality of $\eta_d$ by induction. Since $Dw(0)=0$, the coefficient of $s$ in the exponential is $-kx$. Thus $\gamma_1$ has the stated form. Its denominator has a simple zero at $k=A$, while $-Ax\neq 0$.
\end{proof}

Expand the pointwise factor as
\[
q_{\varphi_x(s)}=e^{-k\varphi_x(s)}\frac{\det g_0}{\det g_{\varphi_x(s)}}=1+\sum_{i\geq1}s^i r_i(k).
\]
Each $r_i(k)$ is polynomial in $k$.
Define
\[
S_d(k)=[s^d]\operatorname{pr}_{d\Lambda_q}\bigl(q_{\varphi_x(s)}\sigma_m(G(s)^{-1})\bigr).
\]
Set $\eta_0=1$ and $A_{\mathbf0}=1$. The highest-weight bounds in Lemma~\ref{lem:highest-weight-support}, together with \eqref{eq:highest-product-projection} and \eqref{eq:pre-fixed-symbol}, give
\begin{equation}\label{eq:projected-density-recursion}
S_d=A_{d(1^q)}(k-p)\eta_d+\sum_{i=1}^{d}\bigl(\operatorname{pr}_{i\Lambda_q}r_i\bigr)\odot\bigl(A_{(d-i)(1^q)}(k-p)\eta_{d-i}\bigr).
\end{equation}
Indeed, in a product of coefficients of degrees $i$ and $d-i$,
only their projections onto $i\Lambda_q$ and $(d-i)\Lambda_q$
can contribute to $d\Lambda_q$.

We shall use the rational function
\begin{equation}\label{eq:rational-density-coefficient}
R_{q,d}(k;x)\coloneq d_{k-p}S_d(k).
\end{equation}
For sufficiently large positive integers $k$, \eqref{eq:varyingbergman} gives
\begin{equation}\label{eq:rational-density-identification}
R_{q,d}(k;x)=[s^d]\operatorname{pr}_{d\Lambda_q}\mathcal B_k^{K_M}(h_{\varphi_x(s)};\cdot).
\end{equation}
For complex $x$, the right-hand side is understood via the complex multilinear extension of the Taylor coefficients at $h_0$.

We now show that one of these rational functions must have a nonzero pole at $k=A$, and the required Taylor order depends on $A$. Set
\[
d_q=
\begin{cases}
2,& A\geq 2,\\
q,& A=1.
\end{cases}
\]

\begin{proposition}\label{prop:nonzero-density-pole}
For every $2\leq q\leq r$ and every $0\neq x\in\mathcal V_{(1^q)}$, at least one of the rational functions
\[
R_{q,1}(k;x),\ldots,R_{q,d_q}(k;x)
\]
has a pole at $k=A$.
\end{proposition}

\begin{proof}
Suppose that all these functions are holomorphic at $k=A$. Since $d_{k-p}$ has a simple zero there,
\[
\operatorname{pole}_{k=A}S_d\leq 1,
\qquad 1\leq d\leq d_q.
\]
In \eqref{eq:projected-density-recursion}, every $r_i$ and every $A_{(d-i)(1^q)}(k-p)$ is holomorphic at $A$, and $A_{d(1^q)}(A-p)\neq 0$.
Induction on $d$ therefore gives
\begin{equation}\label{eq:eta-simple-poles}
\operatorname{pole}_{k=A}\eta_d\leq1,
\qquad 1\leq d\leq d_q.
\end{equation}
Let
\[
g=[\gamma_1]_{-1}\neq 0.
\]
By Lemma~\ref{lem:pre-cartan-nonvanishing},
$g^{\odot j}\neq 0$ for every $j\geq 1$.

If $A\geq 2$, then $\nu_1(q)=\nu_2(q)=1$, whereas
\[
\eta_2=-\gamma_2+\gamma_1^{\odot2}.
\]
Consequently,
\[
[\eta_2]_{-2}=g^{\odot2}\neq 0,
\]
contrary to \eqref{eq:eta-simple-poles}.

Now suppose that $A=1$. Formula \eqref{eq:denominator-zero-order} gives
\[
\nu_j(q)=j\quad(1\leq j<q),\qquad \nu_q(q)=q-1.
\]
We claim that
\begin{equation}\label{eq:gamma-leading-coefficients}
[\gamma_j]_{-j}=g^{\odot j},\qquad 1\leq j<q.
\end{equation}
The case $j=1$ is the definition of $g$. Suppose $2\leq j<q$ and the claim holds through degree $j-1$. For $i\leq j-2$, the product $\gamma_i\odot\eta_{j-i}$ has pole order at most $i+1\leq j-1$. Since $\eta_1=-\gamma_1$ and $\eta_j$ has at most a simple pole, the coefficient of $(k-A)^{-j}$ in \eqref{eq:eta-recursion} gives
\[
0=-[\gamma_j]_{-j}+[\gamma_{j-1}]_{-(j-1)}\odot g.
\]
This proves \eqref{eq:gamma-leading-coefficients}.

At degree $q$, the term $\gamma_q$ has pole order at most $q-1$. All products in \eqref{eq:eta-recursion}, except $\gamma_{q-1}\odot\eta_1$, also have pole order at most $q-1$.
Thus
\[
[\eta_q]_{-q}=[\gamma_{q-1}]_{-(q-1)}\odot g=g^{\odot q}\neq 0,
\]
again contradicting \eqref{eq:eta-simple-poles}.
\end{proof}

\begin{remark}
Proposition~\ref{prop:nonzero-density-pole} already rules out nonconstant real-analytic curves of Bergman logarithmically flat metrics on the reduced graph. Indeed, logarithmic flatness forces every rational function $R_{q,d}(k;x)$ arising along such a curve to have the form $P(k)+c/k$, contradicting the nonzero pole at $k=a/2$ established above. In the next section, we strengthen this argument to obtain local rigidity from the vanishing of finitely many Bergman coefficients.
\end{remark}

\section{Local rigidity when $a$ is even}
\label{sec:local-rigidity-even}

We now complete the proof of Theorem~\ref{thm:localrigidity} when $a$ is even. We assume that $r\geq 2$.

\begin{lemma}\label{lem:common-denominator}
Define
\begin{equation}\label{eq:common-denominator}
Q(k)=\prod_{q=2}^{r}\prod_{j=1}^{d_q}D_{j\mathbf (1^q)}(k-p)^{\lfloor d_q/j\rfloor}.
\end{equation}
For every $2\leq q\leq r$, $1\leq d\leq d_q$, and $x\in\mathcal V_{(1^q)}$, the function $Q(k)R_{q,d}(k;x)$ is a polynomial in $k$.
\end{lemma}

\begin{proof}
The numerator in \eqref{eq:gamma-explicit} is a polynomial in $k$. Iterating \eqref{eq:eta-recursion} gives
\[
\eta_h=\sum_{l=1}^{h}(-1)^l\sum_{\substack{j_1+\cdots+j_l=h\\j_i\geq1}}\gamma_{j_1}\odot\cdots\odot\gamma_{j_l}.
\]
In a term with $h\leq d_q$, the index $j$ occurs at most $\lfloor d_q/j\rfloor$ times. Hence the factor of $Q$ indexed by $q$ clears the denominators of all $\eta_h$, $1\leq h\leq d_q$. All remaining factors in \eqref{eq:projected-density-recursion} and \eqref{eq:rational-density-coefficient} are polynomials in $k$.
\end{proof}

Let $N_0$ be the cutoff chosen in the elliptic reduction and set
\begin{equation}\label{eq:finite-coefficient-cutoff}
L_Q=\max\{2,\deg Q+1,N_0-n\},\qquad J=n+L_Q.
\end{equation}
In particular, $J\geq N_0$.

\begin{lemma}\label{lem:finite-laurent-vanishing}
Let $R(k)$ be a rational function with values in a finite-dimensional complex vector space, and suppose that $Q(k)R(k)$ is a polynomial. If
\[
[k^{-l}]R(k)=0,\qquad 2\leq l\leq L_Q,
\]
then
\[
R(k)=P(k)+\frac{c}{k}
\]
for a vector-valued polynomial $P$ and a constant vector $c$. In particular, $R$ is holomorphic at $k=A$.
\end{lemma}

\begin{proof}
Let $P$ be the polynomial part of the Laurent expansion of $R$ at infinity, and let $c=[k^{-1}]R$. The remainder $R-P-c/k$ has denominator dividing $kQ(k)$ and satisfies
\[
R(k)-P(k)-\frac{c}{k}=O(k^{-L_Q-1}).
\]
Multiplication by $kQ(k)$ gives a polynomial tending to zero at infinity, because $L_Q\geq\deg Q+1$. That polynomial is identically zero.
\end{proof}

Fix an integer
\begin{equation}\label{eq:even-regularity}
s\geq\max\{4N_0+4,\,2J+2\},
\end{equation}
and choose the compatible graph map $w$ with values in $C^{s,\alpha}(M,\mathbb R)$, as constructed in Section~\ref{sec:elliptic-reduction}.
On a sufficiently small neighborhood $U$ of $0$ in $E$, define
\[
\Phi(v)=\left(\operatorname{pr}_{d\Lambda_q}\mathcal A_j(v+w(v))\right)_{\substack{2\leq q\leq r,\ 1\leq d\leq d_q\\n+2\leq j\leq J}}.
\]
This is a real-analytic map with values in the finite-dimensional real vector space
\[
\bigoplus_{q=2}^r\ \bigoplus_{d=1}^{d_q}\bigoplus_{j=n+2}^{J}\mathcal V_{d(1^q)}^{\mathbb R}.
\]

We use the Bergman expansion with derivatives in a smooth real parameter. For every smooth family $\widetilde\varphi(t)$ of potentials defining positive metrics, and every fixed $N\geq 0$,
\begin{equation}\label{eq:finite-order-bergman-expansion}
[t^N]\mathcal B_k^{K_M}(h_{\widetilde\varphi(t)};\cdot)\sim\sum_{j\geq 0}k^{n-j}[t^N]\mathcal A_j(\widetilde\varphi(t)).
\end{equation}
Here $[t^N]$ means the Taylor coefficient at $t=0$. This follows from the parameter estimates in \cite[Theorems~1.1 and~1.3]{DaiLiuMa2006}.

\begin{proposition}\label{prop:isolated-projected-zero}
Zero is an isolated zero of $\Phi$.
\end{proposition}

\begin{proof}
Suppose otherwise. Since $\Phi$ is real analytic and the origin is not an isolated point of $\Phi^{-1}(0)$, the real-analytic curve selection lemma (see, e.g., \cite[p.~103]{Coltoiu1997}) yields a real-analytic curve $v\colon (-\varepsilon,\varepsilon)\to E$ such that
\[
v(0)=0,\qquad v(t)\neq 0,\qquad \Phi(v(t))=0\quad (0<t<\varepsilon).
\]
Choose $q$, $e_q$, and $x_q\neq 0$ as in Lemma~\ref{lem:leading-highest-weight-coefficient}. Put $\varphi(t)=v(t)+w(v(t))$. For $1\leq d\leq d_q$, set $N=de_q$ and let
\[
\widetilde\varphi_N(t)=\sum_{l=1}^{N}\frac{t^l}{l!}\varphi^{(l)}(0).
\]
By the smoothness of the Taylor tensors of $w$, the truncated curve $\widetilde\varphi_N(t)$ consists of smooth potentials and defines positive K\"ahler metrics for small $|t|$. Its Taylor expansion agrees with that of $\varphi(t)$ through order $N$.

For sufficiently large integers $k$, Lemma~\ref{lem:leading-highest-weight-coefficient} and \eqref{eq:rational-density-identification} give
\[
[t^N]\operatorname{pr}_{d\Lambda_q}\mathcal B_k^{K_M}(h_{\widetilde\varphi_N(t)};\cdot)=R_{q,d}(k;x_q).
\]
Applying \eqref{eq:finite-order-bergman-expansion} to the left-hand side and comparing Laurent coefficients at infinity yields
\[
[k^{-l}]R_{q,d}(k;x_q)=[t^N]\operatorname{pr}_{d\Lambda_q}\mathcal A_{n+l}(\varphi(t)),\qquad 2\leq l \leq L_Q.
\]
In this identity, replacing $\widetilde\varphi_N$ by $\varphi$ does not change the coefficient: their derivatives through order $N$ agree, and $\mathcal A_{n+l}$ is analytic on the chosen H\"older space. Since $\Phi(v(t))=0$, all the displayed coefficients vanish. Lemmas~\ref{lem:common-denominator} and \ref{lem:finite-laurent-vanishing} therefore imply that $R_{q,d}(k;x_q)$ is holomorphic at $k=A$ for every $1\leq d\leq d_q$. This contradicts Proposition~\ref{prop:nonzero-density-pole}.
\end{proof}

We are now in a position to prove Theorem~\ref{thm:localrigidity} when $a$ is even. We restate this case below for the reader's convenience.

\begin{theorem}\label{thm:even-a-rigidity}
Let $M$ be an irreducible Hermitian symmetric space of compact type of complex dimension $n$, equipped with its standard polarization $(L,h_0)$. Assume that the restricted-root parameter $a$ is even. For every $\alpha\in(0,1)$, there exist integers $J\geq n+2$ and $s\geq 2$, and a constant $\varepsilon>0$, such that the following holds.

If $\varphi\in C^\infty(M,\mathbb R)$ satisfies
\[
h_\varphi=h_0e^{-\varphi},\qquad\omega_\varphi>0,\qquad\|\varphi\|_{C^{s,\alpha}(M,\omega_0)}<\varepsilon,
\]
and
\[
a^{K_M}_j(h_\varphi)\equiv0\qquad(n+2\leq j\leq J),
\]
then
\[
h_\varphi=CF^*h_0
\]
for some $C>0$ and $F\in\operatorname{Aut}^0(M,L)$. In particular, the conclusion holds if the corresponding unit circle bundle is Bergman logarithmically flat.
\end{theorem}
\begin{proof}
For $r\geq 2$, choose $J$ and $s$ as in \eqref{eq:finite-coefficient-cutoff} and \eqref{eq:even-regularity}. Let $\varphi$ satisfy the hypotheses of the theorem with these choices. By Lemma~\ref{lem:local-moment-slice}, choose a polarized
automorphism $F$ and a real constant $c$ so that
\[
h_\psi=e^{-c}F^*h_\varphi,\qquad \psi\perp\mathbb R\mathbf1\oplus\mathcal V_{(1)}^{\mathbb R}.
\]
For sufficiently small $\|\varphi\|_{C^{s,\alpha}}$, the normalized potential $\psi$ is also small in $C^{s,\alpha}$. The coefficient equations are preserved by this transformation,
so
\[
\mathcal A_j(\psi)=0,\qquad n+2\leq j\leq J.
\]
Since $J\geq N_0$, the elliptic reduction gives
\[
\psi=v+w(v),\qquad \text{for some }v\in E.
\]
The coefficient equations imply $\Phi(v)=0$. Proposition~\ref{prop:isolated-projected-zero} gives $v=0$, and hence $\psi=0$. It follows that
\[
h_\varphi=e^c(F^{-1})^*h_0.
\]
This proves the required conclusion. The assertion for Bergman logarithmic flatness follows because it implies $a^{K_M}_j(h_\varphi)=0$ for every $j\geq n+2$.
\end{proof}

\section{Local rigidity from later Bergman coefficients}
\label{sec:later-bergman-coefficients}

We prove Theorem~\ref{thm:finite-block-rigidity} by the argument of Sections~\ref{sec:infinitesimal-deformations}--\ref{sec:local-rigidity-even}, indicating only the changes needed when the coefficient block starts at a fixed $q\geq n+2$.

First, the common linearized kernel remains unchanged. Indeed, a rational function has vanishing Laurent coefficients of $k^{n-j}$ for all $j\geq q$ if and only if it has the form
\begin{equation}\label{eq:later-laurent-form}
 R(k)=P(k)+\sum_{\nu=1}^{q-n-1}c_\nu k^{-\nu},
\end{equation}
where $P$ is polynomial. Thus the criterion in Proposition~\ref{prop:divisibility-criterion} becomes
\[
D_{\mathbf s}(m)\mid (m+p)^{q-n}d_m.
\]
The additional power of $m+p$ does not affect the nonzero poles used in the proof of Theorem~\ref{thm:linearization-repeat}.

Since $\mathcal P_q$ is elliptic by Lemma~\ref{lem:coefficient-analyticity}, choose $N_q\geq q$ such that
\[
 \bigcap_{j=q}^{N_q}\ker\mathcal P_j=\mathcal K_{\operatorname{lin}}.
\]
The reduction in Section~\ref{sec:elliptic-reduction} applies to $\mathcal F_q=\sum_{j=q}^{N_q}\mathcal P_j\mathcal A_j$ and gives
an equivariant analytic graph $v+w_q(v)$ over the same residual space $E$. If $E=0$, the proof is complete. Otherwise, the argument
of Section~\ref{sec:nonlinear-obstructions} applies without change: for each nonzero residual direction, one of the selected rational
density coefficients has a pole at $k=a/2$.

It remains to modify the finite Laurent test in Lemma~\ref{lem:finite-laurent-vanishing}. If $QR$ is polynomial and
\[
 [k^{-\nu}]R=0\qquad(q-n\leq\nu\leq J-n), \qquad J\geq q+\deg Q-1,
\]
then $R$ has the form \eqref{eq:later-laurent-form}. To see this, subtract its polynomial part and its first $q-n-1$ negative-power terms. Multiplication of the remainder by $k^{q-n-1}Q(k)$ gives a polynomial tending to zero at infinity, hence identically zero.

Take $Q$ as in \eqref{eq:common-denominator}, and choose
\[
 J=\max\{N_q,q,q+\deg Q-1\},\qquad s\geq\max\{4N_q+4,2J+2\}.
\]
The proof of Proposition~\ref{prop:isolated-projected-zero}, using the coefficient block $\mathcal A_q,\ldots,\mathcal A_J$ and the
preceding Laurent test, again excludes every nonconstant analytic arc in the reduced zero set: \eqref{eq:later-laurent-form} is
regular at $k=a/2$, contradicting Proposition~\ref{prop:nonzero-density-pole}. Applying Lemma~\ref{lem:local-moment-slice} then we complete the proof of Theorem~\ref{thm:finite-block-rigidity}.

\section{Products and flat torus factors}
\label{sec:products}

We now prove Theorem~\ref{thm:generalrigidity} using the notation introduced in \eqref{eq:general-product-polarization} and \eqref{eq:general-reference-metric}. The proof follows the strategy we developed above for irreducible compact Hermitian symmetric spaces in Sections~\ref{sec:infinitesimal-deformations}--\ref{sec:local-rigidity-even}. We first determine the common kernel of the linearized Bergman coefficient equations for the product polarization, and then adapt the finite-dimensional elliptic reduction and the nonlinear obstruction argument to the product setting. We give the details of these modifications, while referring to the earlier sections for the parts of the argument that remain unchanged.

\subsection{The product linearization}

First let $X=\prod_{i=1}^{l}M_i$ be of compact type, with polarization $L_X=L_1^{d_1}\boxtimes\cdots\boxtimes L_l^{d_l}$ and reference invariant metric $h_X$. Put
\[
 N=\dim_{\mathbb C}X,\qquad R=\sum_i r_i,\qquad m_i=d_i k-p_i,\qquad c_{ij}=\frac{a_i}{2}(j-1),
\]
and fix $q\geq N+2$. For the product labels $\mathbf s=(\mathbf s^{(1)},\ldots,\mathbf s^{(l)})$ such that $\mathbf s^{(i)}\in \mathcal P_{r_i}$ for each $i$, the product Peter--Weyl decomposition is
\[
 L^2(X,\mathbb C)=\widehat{\bigoplus}_{\mathbf s}\mathcal V_{\mathbf s}, \qquad \mathcal V_{\mathbf s} =\bigotimes_i\mathcal V_{\mathbf s^{(i)}}^{(i)}.
\]
Define
\[
 A_{\mathbf s}(k)=\prod_i A_{\mathbf s^{(i)}}^{(i)}(m_i), \qquad D_{\mathbf s}(k)=\prod_i D_{\mathbf s^{(i)}}^{(i)}(m_i),
\]
using the factor polynomials in \eqref{eq:numerator} and \eqref{eq:denominator}. The Hilbert polynomial is
\[
 H(k)=C_X\prod_{i=1}^{l}\prod_{j=1}^{r_i} \prod_{u=0}^{b_i+a_i(r_i-j)} \bigl(d_i k-(1+c_{ij}+u)\bigr),\qquad C_X>0.
\]
Let $T_{X,k}$ denote the Toeplitz map associated with $H^0(X,L_X^k\otimes K_X)$ at the reference metric. By Proposition~\ref{prop:varying-identity}, the Gram matrix with respect to $h_{\varphi}$ is $G_\varphi=T_{X,k}(e^{-k\varphi})$. Thus the first variational formula in Theorem~\ref{thm:first-variation} gives the scalar
\begin{equation}\label{eq:prod-linear-scalar}
 H(k)\left(k\frac{A_{\mathbf s}(k)}{D_{\mathbf s}(k)}-k+\lambda_{\mathbf s}\right),\qquad\lambda_{\mathbf s}=\sum_id_i^{-1}\lambda_{\mathbf s^{(i)}}^{(i)}.
\end{equation}
We use $\mathcal A_j$ and $\mathcal P_j$ as in \eqref{eq:coefficient-equation}, and write $\mathcal K_q(X)=\bigcap_{j\geq q}\ker\mathcal P_j$.

For $\mathbf h=(h_1,\ldots,h_l)$ with $0\leq h_i\leq r_i$, put
\[
 (1^{\mathbf h})=((1^{h_1}),\ldots,(1^{h_l})),\qquad u(\mathbf h)=\#\{i\colon h_i>0\},
\]
where $(1^0)=\mathbf0$, and define
\[
 D_{\mathbf h}^{+}(k)=\prod_{i=1}^{l}\prod_{j=2}^{h_i}\left(k-\frac{c_{ij}}{d_i}\right).
\]

The common kernel of the linearized Bergman coefficient equations here is the following.
\begin{proposition}\label{prop:prod-linear-kernel}
For every $q\geq N+2$,
\[
 \mathcal K_q(X)=\bigoplus_{\substack{0\leq h_i\leq r_i,\ u(\mathbf h)\leq q-N\\D_{\mathbf h}^{+}\mid H A_{(1^{\mathbf h})}}}\mathcal V_{(1^{\mathbf h})}^{\mathbb R}.
\]
The summands with $\sum_i h_i\leq 1$ form the space of infinitesimal symmetry directions.
\[
\mathcal S_X=\mathbb R\mathbf1\oplus \bigoplus_i\pi_i^*\mathcal V_{(1),i}^{\mathbb R}.
\]
The remaining summands form $E=\mathcal K_q(X)\cap\mathcal S_X^\perp$.
\end{proposition}

\begin{proof}
By \eqref{eq:later-laurent-form} and \eqref{eq:prod-linear-scalar}, $\mathcal V_{\mathbf s}$ belongs to the complexified kernel precisely when 
\begin{equation}\label{eq:prod-linear-divisibility}
D_{\mathbf s}(k)\ \bigm|\ k^{q-N}H(k)A_{\mathbf s}(k).
\end{equation}
All roots of $H$ and $A_{\mathbf s}$ are positive. If $s_1^{(i)}\geq2$, the denominator has the uncancelled root $-1/d_i$. Thus every admissible label has the form ${(1^\mathbf h})=((1^{h_1}),\ldots,(1^{h_l}))$, where $0\leq h_i\leq r_i$. For such a tuple,
\[
 D_{(1^{\mathbf h})}(k)=\left(\prod_i d_i^{h_i}\right)k^{u(\mathbf h)}D_{\mathbf h}^{+}(k).
\]
Since $H(0)A_{(1^{\mathbf h})}(0)\neq 0$, the asserted conditions are equivalent to \eqref{eq:prod-linear-divisibility}. The symmetry statement follows from Proposition~\ref{prop:constants-holomorphy}.
\end{proof}

\subsection{Local deformations}

Since $\Delta_X=\sum_i d_i^{-1}\Delta_i$ and $m_i=d_i k-p_i$, the coefficient of $k^{-j}$ in the product Berezin expansion has principal differential part
\[
\sum_{a_1+\cdots+a_l=j}\prod_{i=1}^{l}\frac{(-d_i^{-1}\Delta_i)^{a_i}}{a_i!}=\frac{(-\Delta_X)^j}{j!}.
\]
The analyticity, ellipticity, and self-adjointness arguments of Lemma~\ref{lem:coefficient-analyticity} therefore can apply. Choose $N_q\geq q$ so that the common linearized kernel has stabilized. The argument in Section~\ref{sec:elliptic-reduction} gives an equivariant graph
\[
\varphi=v+w(v),\qquad v\in E,\qquad w(0)=Dw(0)=0,
\]
with smooth Taylor tensors at $0$. Let $\mathcal I$ be the set of tuples $\mathbf h=(h_1,\cdots,h_l)$ indexing the summands of
$E^{\mathbb C}$, so that
\[
E^{\mathbb C}=\bigoplus_{\mathbf h\in\mathcal I}\mathcal V_{(1^{\mathbf h})},\qquad (1^{\mathbf h})=((1^{h_1}),\ldots,(1^{h_l})).
\]
If $E=0$, this proves the compact type assertion. Otherwise, we need the following observation in the new pole argument.

\begin{lemma}\label{lem:prod-inverse-series}
Let $A=\bigoplus_{d\geq0}A_d$ be a graded commutative complex integral domain such that $A_0=\mathbb C$, and each $A_d$ is finite-dimensional. For each $d\geq1$, let $\gamma_d(\zeta)$ and $\eta_d(\zeta)$ be $A_d$-valued functions which are meromorphic near $\zeta=0$. Suppose that
\[
\left(1+\sum_{d\geq1}\gamma_d(\zeta)z^d\right)\left(1+\sum_{d\geq1}\eta_d(\zeta)z^d\right)=1
\]
as formal power series in $z$ and $\gamma_1$ has a pole at $\zeta=0$. Then the following holds: for integers $B\geq1$ and $H_0\geq0$, the bounds
\[
\operatorname{pole}_0\gamma_d\leq B,\qquad \operatorname{pole}_0\eta_d\leq H_0
\qquad (1\leq d\leq B+H_0)
\]
cannot all hold.
\end{lemma}

\begin{proof}
Suppose that all these bounds hold, and put $M=B+H_0$. Set
\[
\theta=\max_{1\leq d\leq M}\left\{\frac{\operatorname{pole}_0\gamma_d}{d},\frac{\operatorname{pole}_0\eta_d}{d}\right\}.
\]
Since $\gamma_1$ has a pole at $0$, $\theta\geq1$. By definition, the Laurent expansions of $\gamma_d$ and $\eta_d$ contain no powers of $\zeta$ with exponent less than $-\theta d$.

For $1\leq d\leq M$, define
\[
a_d=
\begin{cases}
[\zeta^{-\theta d}]\gamma_d(\zeta),&\theta d\in\mathbb Z,\\
0,&\theta d\notin\mathbb Z,
\end{cases}
\qquad
b_d=
\begin{cases}
[\zeta^{-\theta d}]\eta_d(\zeta),&\theta d\in\mathbb Z,\\
0,&\theta d\notin\mathbb Z.
\end{cases}
\]
Put $a_0=b_0=1$ and
\[
g(z)=\sum_{d=0}^{M}a_dz^d,\qquad h(z)=\sum_{d=0}^{M}b_dz^d.
\]
By the definition, at least one of $g,h$ is nonconstant.

We claim that
\[
\sum_{i=0}^{n}a_i b_{n-i}=0\qquad (1\leq n\leq M).
\]
Set $\gamma_0=\eta_0=1$. The assumed inverse relation gives
\[
\sum_{i=0}^{n}\gamma_i(\zeta)\eta_{n-i}(\zeta)=0.
\]
If $\theta n$ is an integer, compare the coefficients of $\zeta^{-\theta n}$ in this identity. Every Laurent exponent in $\gamma_i$ is at least $-\theta i$, and every Laurent exponent in $\eta_{n-i}$ is at least $-\theta(n-i)$. Their sum can equal $-\theta n$ only when both inequalities
are equalities. Consequently,
\[
[\zeta^{-\theta n}]\bigl(\gamma_i(\zeta)\eta_{n-i}(\zeta)\bigr)=a_i b_{n-i},
\]
which proves the required identity. If $\theta n$ is not an integer, each product $a_i b_{n-i}$ is zero.

It follows that
\[
g(z)h(z)\equiv1\pmod{z^{M+1}}.
\]
On the other hand, if $a_d\neq 0$ for $d\geq1$, then $\gamma_d$ has a pole of order $\theta d$, so $\theta d\leq B$. Similarly, $b_d\ne0$ implies $\theta d\leq H_0$. Hence
\[
\deg g+\deg h\leq\left\lfloor\frac{B}{\theta}\right\rfloor+\left\lfloor\frac{H_0}{\theta}\right\rfloor\leq B+H_0=M.
\]
Thus $gh$ has degree at most $M$, which implies that $gh=1$. Since $A$ is an integral domain,
\[
\deg(gh)=\deg g+\deg h.
\]
Both $g$ and $h$ must therefore be constant, contradicting the construction.
\end{proof}

For each $\mathbf h\in\mathcal I$, consider the graded vector space
\[
A_{(1^{\mathbf h})}=\bigoplus_{d\geq0}\mathcal V_{d(1^{\mathbf h})},
\]
equipped with the Cartan product
\[
f\odot g=\operatorname{pr}_{\mathcal V_{(d+e)(1^{\mathbf h})}}(fg),\qquad f\in\mathcal V_{d(1^{\mathbf h})},\quad g\in\mathcal V_{e(1^{\mathbf h})}.
\]
Then $A_{(1^{\mathbf h})}$ is a graded commutative algebra.

\begin{lemma}\label{lem:prod-cartan-domain}
The algebra $A_{(1^{\mathbf h})}$ is an integral domain.
\end{lemma}

\begin{proof}
Fix nonzero elements $v\in\mathcal V_{d(1^{\mathbf h})}$ and $w\in\mathcal V_{e(1^{\mathbf h})}$. Choose nonzero highest-weight
covectors
\[
\ell_d\in\mathcal V_{d(1^{\mathbf h})}^*,\qquad \ell_e\in\mathcal V_{e(1^{\mathbf h})}^*,\qquad \ell_{d+e}\in\mathcal V_{(d+e)(1^{\mathbf h})}^*.
\]
Let $P_{d,e}\colon\mathcal V_{d(1^{\mathbf h})}\otimes\mathcal V_{e(1^{\mathbf h})}\longrightarrow\mathcal V_{(d+e)(1^{\mathbf h})}$ be the Cartan product map, which is nonzero. Then $P_{d,e}^*\ell_{d+e}$ and $\ell_d\otimes\ell_e$ span the same highest-weight line. After rescaling $\ell_{d+e}$, we have
\[
\ell_{d+e}(u\odot u')=\ell_d(u)\ell_e(u')
\]
for all $u\in\mathcal V_{d(1^{\mathbf h})}$ and $u'\in\mathcal V_{e(1^{\mathbf h})}$.

Let $G^*=\prod_{i=1}^{l}G_i^*$ be the connected compact group acting on $X=\prod_{i=1}^{l}M_i$ by the product action. The functions
\[
F_v(g)=\ell_d(g\cdot v),\qquad F_w(g)=\ell_e(g\cdot w),\qquad g\in G^*,
\]
are real analytic and not identically zero. And their product is not identically zero. By $G^*$-equivariance,
\[
\ell_{d+e}\bigl(g\cdot(v\odot w)\bigr)=F_v(g)F_w(g),
\]
so $v\odot w\neq 0$.

Now let $F,G\in A_{(1^{\mathbf h})}$ be nonzero, with highest-degree components $f_D$ and $g_E$. The component of degree $D+E$ in $F\odot G$ is $f_D\odot g_E\neq 0$. Hence $F\odot G\neq 0$.
\end{proof}

Let $\Lambda_{(1^{\mathbf h})}$ be the highest weight of $\mathcal V_{(1^{\mathbf h})}$. Fix $0\neq x\in\mathcal V_{(1^{\mathbf h})}$ and put $\varphi_x(z)=zx+w(zx)$, using the complexified Taylor tensors of $w$. Denote by $\gamma_d,\eta_d,R_{\mathbf h,d}$ the product versions of \eqref{eq:gamma-explicit}, \eqref{eq:eta-recursion}, and \eqref{eq:rational-density-coefficient}. The same argument gives these identities. In particular, $\gamma_d$ has denominator $D_{d(1^{\mathbf h})}(k)$ and polynomial numerator divisible by $k$, and
\[
 \gamma_1(k)=-\frac{kx}{D_{(1^{\mathbf h})}(k)}.
\]
The Bergman kernel coefficients satisfy the identity \eqref{eq:projected-density-recursion}, with the product polynomials and with $H(k)$ replacing $d_{k-p}$.

Put $D_*=R+q-N-2$. We claim that, for some $1\leq d\leq D_*$, $R_{\mathbf h,d}$ does not have the form
\begin{equation}\label{eq:prod-allowed-laurent}
P(k)+\sum_{\nu=1}^{q-N-1}c_\nu k^{-\nu}.
\end{equation}

Suppose first that $\mathbf h$ has exactly one nonzero entry, with $h_i\geq 2$. Put $k_0=a_i/(2d_i)>0$. The factor formulae give $A_{d(1^{\mathbf h})}(k_0)\neq 0$, since $p_i-c_{i h_i}>a_i/2$, and $\gamma_1$ has a simple pole there. The denominator factors with index $j=1$ do not vanish at $k_0$, while those with each index $2\leq j\leq h_i$ contribute
at most one zero. Hence
\[
\operatorname{pole}_{k_0}\gamma_d\leq h_i-1.
\]
Put $H_0=\operatorname{ord}_{k_0}H$. If all Bergman kernel coefficients of degrees $1\leq d\leq h_i-1+H_0$ had the form \eqref{eq:prod-allowed-laurent}, they would be holomorphic at $k_0$. Then \eqref{eq:projected-density-recursion} would then imply
\[
\operatorname{pole}_{k_0}\eta_d\leq H_0\qquad (1\leq d\leq h_i-1+H_0),
\]
contradicting Lemma~\ref{lem:prod-inverse-series} with $B=h_i-1$.

In the Hilbert polynomial of the $i$-th factor, only the factors with index $j=1$ can vanish at $k_0$, contributing at most one zero. Each other $\nu$-th factor contributes at most $r_\nu$ zeros, counted with multiplicity. Thus
\[
H_0\leq1+\sum_{\nu\neq i}r_\nu,
\]
and consequently
\[
h_i-1+H_0\leq h_i+\sum_{\nu\neq i}r_\nu\leq R\leq D_*.
\]
Suppose next that $\mathbf h$ has at least two nonzero entries. At $k=0$, each product $\prod_{v=0}^{d-1}(d_i k-c_{ij}+v)$ in $D_{d(1^{\mathbf h})}$ has at most a simple zero. Since the numerator of $\gamma_d$ is divisible by $k$,
\[
\operatorname{pole}_0\gamma_d\leq\sum_i h_i-1\leq R-1.
\] 
Thus
\[
\operatorname{pole}_0\gamma_d\leq R-1,\qquad \operatorname{pole}_0\gamma_1=u(\mathbf h)-1>0.
\]
Since $H(0)A_{d(1^{\mathbf h})}(0)\neq 0$, the form \eqref{eq:prod-allowed-laurent} for all $d\leq D_*$ would imply $\operatorname{pole}_0\eta_d\leq q-N-1$ by \eqref{eq:projected-density-recursion}. Lemma~\ref{lem:prod-inverse-series}, with $B=R-1$ and $H_0=q-N-1$, gives a contradiction by degree $D_*$.

The remaining argument is the analogue of Proposition~\ref{prop:isolated-projected-zero}. The weights $\Lambda_{(1^{\mathbf h})}$ are distinct sums, with coefficients zero or one, of the linearly independent strongly orthogonal roots of the factors. Consequently, no two are positive scalar multiples, so Lemma~\ref{lem:leading-highest-weight-coefficient} applies unchanged. Choose a common denominator $Q_X$, independent of $x$, for the selected $R_{(1^{\mathbf h}),d}$ as in Lemma~\ref{lem:common-denominator}, and take
\[
 J_X=\max\{N_q,q,q+\deg Q_X-1\}.
\]
At sufficiently high H\"older regularity, apply the same isolation argument
to
\[
 \operatorname{pr}_{d\Lambda_{(1^{\mathbf h})}}\mathcal A_j(v+w(v)),\qquad (1^{\mathbf h})\in\mathcal I,\quad 1\leq d\leq D_*, \quad q\leq j\leq J_X.
\]
The finite Laurent test in Section~\ref{sec:later-bergman-coefficients} would force \eqref{eq:prod-allowed-laurent} along any nonconstant analytic arc, contrary to the claim. Applying Lemma~\ref{lem:local-moment-slice} completes the proof when $X$ is of compact type.

\subsection{Flat torus factors}

Now let $M=T\times X$, with $\dim_{\mathbb C}T=t$, $n=t+N$, and $q\geq n+2$. In flat holomorphic coordinates and a suitable local
holomorphic frame $e$, the metric satisfies $h_T(e,e)=e^{-|z|^2}$. By locality of the Berezin expansion, its coefficients therefore agree with those of the Bargmann--Fock model; see \cite[\S2.8.1]{Englis2006Quantization}. Hence
\begin{equation}\label{eq: fock-bergman}
B_{T,k}\sim\sum_{j\geq 0}\frac{(-\Delta_T)^j}{j!\,k^j}=\exp(-\Delta_T/k).
\end{equation}
Since the invariant K\"ahler metric on $T$ is flat, all higher Bergman coefficients vanish. Hence
\[
\mathcal B_k^{K_T}(h_T)=C_Tk^t+O(k^{-\infty}).
\]
Write $T=\mathbb C^t/\Lambda$. For each real-linear functional $\xi\colon\mathbb C^t\to\mathbb R$ satisfying $\xi(\Lambda)\subset\mathbb Z$, the function
\[
\chi_\xi(x)=e^{2\pi\ii\xi(x)}
\]
is well defined on $T$ and satisfies $\Delta_T\chi_\xi=\lambda_\xi\chi_\xi$. For all $\xi\neq 0$, we have $\lambda_\xi>0$.

We use the product measure $d\mu_M=dV_T\otimes d\mu_X$, where $dV_T=\omega_T^t/t!$ and $d\mu_X$ is the invariant probability measure on $X$. This measure differs from the probability measure on $M$ only by a constant factor. By the first variational formula \eqref{eq:first-variation}, the eigenvalue of $\mathcal P_j^M$ on $\mathbb C\chi_\xi\otimes\mathcal V_{\mathbf s}$ is the coefficient of $k^{n-j}$ in the expansion at infinity of
\begin{equation}\label{eq:torus-linear-scalar}
C_Mk^tH(k)\left(k e^{-\lambda_\xi/k}\frac{A_{\mathbf s}(k)}{D_{\mathbf s}(k)}-k+\lambda_\xi+\lambda_{\mathbf s}
\right).
\end{equation}
If these coefficients vanish for all $j\geq q$, the expression has the form
\[
P(k)+\sum_{\nu=1}^{q-n-1}c_\nu k^{-\nu},
\]
where $P$ is a polynomial. For $\xi\neq 0$, however, it has an essential singularity at $k=0$. Hence only $\xi=0$ can occur in the common linearized kernel.

Thus every element of the common kernel is independent of the torus variable. Setting $\xi=0$ in \eqref{eq:torus-linear-scalar} gives $C_Mk^t$ times \eqref{eq:prod-linear-scalar}. Since $n=N+t$, comparison of coefficients yields
\[
\mathcal P_j^M(\pi_X^*f)=C_M\pi_X^*(\mathcal P_j^Xf).
\]
Consequently,
\[
\bigcap_{j\geq q}\ker\mathcal P_j^M=\pi_X^*\mathcal K_q(X),
\]
where $\pi_X\colon M\to X$ is the projection and $\mathcal K_q(X)$ is given by Proposition~\ref{prop:prod-linear-kernel} with the same starting index $q$.

The elliptic reduction of Section~\ref{sec:elliptic-reduction} applies on $M$, with symmetry space $\pi_X^*\mathcal S_X$. The coefficient maps commute with translations on $T$. As translations fix the common linearized kernel, uniqueness in the implicit-function theorem implies that the reduced graph is translation-invariant. Thus every sufficiently small normalized solution is of the form $\varphi=\pi_X^*\varphi_X$.

For such potentials, the product formula for the Bergman kernel gives
\[
a_{j,M}^{K_M}(h_T\boxtimes(h_Xe^{-\varphi_X}))=C_Ma_{j,X}^{K_X}(h_Xe^{-\varphi_X}).
\]
The compact-type result with the same starting index $q$ therefore applies, after increasing the finite cutoff and regularity if necessary. This gives $h_\varphi=CF^*h_0$, where $F$ acts on the compact factors and lifts to the product polarization.

Finally, if $M=T$, \eqref{eq: fock-bergman} gives
\[
 \mathcal P_q =\frac{C_T(-1)^{q+1}}{(q+1)!}\Delta_T^{q+1}.
\]
Its kernel consists of constants. Restricting to normalized potentials, the map $\mathcal P_q\mathcal A_q$ has invertible derivative $\mathcal P_q^2$ at zero. The implicit-function theorem therefore gives rigidity with $J=q$ and $h_\varphi=Ch_0$. The logarithmic-flatness assertion follows from Theorem~\ref{thm:logflat}.

\section{Further discussion}
\label{sec:further-discussion}

We conclude with a discussion of future directions of research emanating from our results.

\subsection{Local rigidity beyond symmetric spaces and global rigidity}
\label{subsec:global-local-rigidity}

The Berezin transform and the Laplacian in our setting are both $G^*$-equivariant. By the multiplicity-free spherical decomposition and Schur's lemma, each acts by a scalar on every spherical module. In particular,
\[
B_m\Delta=\Delta B_m.
\]
This reduces the linearized coefficient equations to scalar equations on the spherical summands. Engli\v{s} \cite{Englis2006Symmetric} proved that, under a nondegeneracy assumption, commutation of the Berezin transform with the Laplacian forces that the manifold is Hermitian symmetric. For general polarized K\"ahler manifolds, the spherical decomposition used here is no longer available. We nevertheless expect local rigidity to hold without the symmetry assumption.

\begin{conjecture}[Local rigidity]\label{conj:local-rigidity}
Let $(M,\omega_0;L,h_0)$ be a compact polarized K\"ahler manifold such that the associated circle bundle $S_0$ is Bergman logarithmically flat. Then every sufficiently small $\varphi\in C^\infty(M,\mathbb R)$, with respect to the $C^\infty$ topology, satisfying $\omega_\varphi>0$ and with $S_\varphi$ Bergman logarithmically flat, satisfies
\[
h_\varphi=CF^*h_0
\]
for some $C>0$ and $F\in\operatorname{Aut}^0(M,L)$.
\end{conjecture}

Another natural question is whether the global rigidity of Bergman logarithmic flatness holds on compact Hermitian symmetric spaces. This leads to the following conjecture.

\begin{conjecture}[Global rigidity]\label{conj:global-rigidity}
Let $(M,\omega_0;L,h_0)$ be a compact Hermitian symmetric space with a fixed invariant positive polarization. If $h$ is a smooth positive Hermitian metric on $L$ such that $S_h$ is Bergman logarithmically flat, then
\[
h=CF^*h_0
\]
for some $C>0$ and $F\in\operatorname{Aut}(M,L)$.
\end{conjecture}

In particular, on $(\mathbb{CP}^n,\mathcal O(1))$, the conjecture would imply that Bergman logarithmic flatness characterizes the Fubini--Study metric up to holomorphic automorphisms. One approach to Conjecture~\ref{conj:global-rigidity} is to prove that Bergman logarithmic flatness should imply constant scalar curvature. Then the uniqueness of constant scalar curvature K\"ahler metrics modulo holomorphic automorphisms \cite{BermanBerndtsson2017} would then give the result. A deeper understanding of the relations among the higher Bergman coefficients may be needed.

\subsection{Quantitative local rigidity}
\label{subsec:Quantitative local rigidity}

A by-product of our proof is the local quantitative control of the Bergman coefficient solutions. Fix the invariant polarization in Theorem~\ref{thm:generalrigidity} and an integer $q\geq n+2$. Let
\[
\mathcal O=\{\Gamma(A,c)\colon A\in\mathfrak q,\ c\in\mathbb R
\text{ sufficiently small}\}
\]
be the local orbit of the zero potential under constant rescalings and polarized automorphisms, where $\Gamma$ and $\mathfrak q$ are defined in Subsection~\ref{subsec:symmetry-directions}.

\begin{corollary}\label{cor:quantitative-rigidity}
Fix $0<\alpha<1$. There exist integers $J\geq q$ and $s\geq 2J+2$, constants $C,\varepsilon>0$, and $0<\theta\leq1$ with the following property. For every smooth function $\varphi$ which satisfies $\|\varphi\|_{C^{s,\alpha}}<\varepsilon$, we have
\begin{equation}\label{eq:quantitative-rigidity}
\operatorname{dist}_{C^{s,\alpha}}(\varphi,\mathcal O)\leq C\left(\sum_{j=q}^{J}\|\mathcal A_j(\varphi)\|_{C^{s-2j-2,\alpha}}\right)^\theta.
\end{equation}
If $E=0$, one can take $\theta=1$.
\end{corollary}

\begin{proof}
Let $r_s(\varphi)=\sum_{j=q}^{J}\|\mathcal A_j(\varphi)\|_{C^{s-2j-2,\alpha}}$. Choose the finite coefficient map $\Phi$ and the compatible graph $v+w(v)$ as in the proof of local rigidity, increasing $J$ and $s$ if necessary. Set
\[
\mathcal K_q=\mathcal S\oplus E,\qquad
\mathcal L=D\mathcal F(0),\qquad d=\operatorname{ord}\mathcal L.
\]
We have seen that the restriction $\mathcal L\colon C^{s,\alpha}\cap\mathcal K_q^\perp\longrightarrow C^{s-d,\alpha}\cap\mathcal K_q^\perp$ is an isomorphism, and the image of $\mathcal F$ lies in $\mathcal K_q^\perp$.

By the implicit-function theorem, we can solve $\mathcal F(v+u)=r$ uniquely as $u=W(v,r)$ near the origin, where $v\in E$ and $u,r\in\mathcal K_q^\perp$ belong to the respective H\"older spaces above. Since $W(v,0)=w(v)$, and the derivative of $W$ with respect to $r$ is locally bounded, we have
\[
\|W(v,r)-w(v)\|_{C^{s,\alpha}}\leq C\|r\|_{C^{s-d,\alpha}}.
\]
Now let $\psi\in\mathcal S^\perp$ be small and put $v=\Pi_E\psi$. Then $\psi-v\in\mathcal K_q^\perp$, and
\[
\psi-v=W(v,\mathcal F(\psi)).
\]
Consequently,
\begin{equation}\label{eq:quantitative-graph-distance}
\|\psi-v-w(v)\|_{C^{s,\alpha}}\leq C\|\mathcal F(\psi)\|_{C^{s-d,\alpha}}\leq Cr_s(\psi).
\end{equation}
For the last inequality, recall that $\mathcal F=\sum_{j=q}^{J_0}\mathcal P_j\mathcal A_j$ for some $J_0\leq J$, with $\operatorname{ord}\mathcal P_j=2j+2$ and $d=4J_0+4$. Thus
\[
\|\mathcal P_j\mathcal A_j(\psi)\|_{C^{s-d,\alpha}}\leq C\|\mathcal A_j(\psi)\|_{C^{s-2j-2,\alpha}}\qquad(q\leq j\leq J_0).
\]
By definition, $\Phi(v)$ consists of finitely many projections of the coefficients $\mathcal A_j(v+w(v))$, with $q\leq j\leq J$. Local Lipschitz continuity of these coefficient maps gives
\[
\begin{aligned}
|\Phi(v)|
&\leq C\sum_{j=q}^{J}\|\mathcal A_j(v+w(v))\|_{C^{s-2j-2,\alpha}}\\
&\leq Cr_s(\psi)+C\|\psi-v-w(v)\|_{C^{s,\alpha}}\\
&\leq Cr_s(\psi).
\end{aligned}
\]
Since $\Phi$ is real analytic and $0$ is an isolated zero point, the \L{}ojasiewicz inequality \cite[Remark~6.5]{BierstoneMilman1988} implies that
\[
|v|\leq C|\Phi(v)|^\theta\leq Cr_s(\psi)^\theta
\]
for some $0<\theta\leq1$. Since $w(0)=Dw(0)=0$, using the equivalence of norms on $E$, and \eqref{eq:quantitative-graph-distance}, we obtain
\[
\begin{aligned}
\|\psi\|_{C^{s,\alpha}}&\leq \|v\|_{C^{s,\alpha}}+\|w(v)\|_{C^{s,\alpha}}+\|\psi-v-w(v)\|_{C^{s,\alpha}}\\
&\leq C\bigl(|v|+|v|^2+r_s(\psi)\bigr)\leq Cr_s(\psi)^\theta.
\end{aligned}
\]
Here we choose the neighborhood so that $r_s(\psi)\leq1$. If $E=0$, then $v=w(v)=0$, so \eqref{eq:quantitative-graph-distance} gives this estimate with $\theta=1$ directly.

For sufficiently small $\varphi$, Lemma~\ref{lem:local-moment-slice} gives a small pair $(A,c)$ such that
\[
e^{-c}F_A^*h_\varphi=h_\psi,\qquad\Pi_{\mathcal S}\psi=0.
\]
Since $F_A$ is close to the identity, its pullback is uniformly bounded on the relevant H\"older spaces. Using $\mathcal A_j(\psi)=F_A^*\mathcal A_j(\varphi)$, we obtain
\[
r_s(\psi)\leq Cr_s(\varphi).
\]
Define $\gamma\in\mathcal O$ by $h_\gamma=e^c(F_A^{-1})^*h_0$, then $\varphi-\gamma=\psi\circ F_A^{-1}$, and hence
\[
\operatorname{dist}_{C^{s,\alpha}}(\varphi,\mathcal O)\leq \|\varphi-\gamma\|_{C^{s,\alpha}}\leq C\|\psi\|_{C^{s,\alpha}}\leq Cr_s(\varphi)^\theta.
\]
This proves \eqref{eq:quantitative-rigidity}.
\end{proof}

A natural question is to determine the optimal exponent $\theta$. We note that if $E\neq 0$, the optimal exponent cannot exceed $1/2$: along $\varphi=tv$ with some $0\neq v\in E$, the Bergman coefficient error is $O(t^2)$, whereas the distance to $\mathcal O$ is bounded below by a positive multiple of $|t|$. To prove that $\theta=1/2$ is sharp, it would suffice to show that the quadratic term of $\Phi$ does not vanish on any nonzero vector in $E$.

\appendix
\section{An example of Bergman logarithmic flatness without Szeg\H{o} logarithmic flatness}
\label{app:blowup}

In this appendix, we show that, for every $n\geq 2$, the unit circle bundle associated with the product Fubini--Study polarization on $\operatorname{Bl}_p\mathbb{CP}^n$ is Bergman logarithmically flat but not Szeg\H{o} logarithmically flat with respect to the natural measure. We note that the Bergman kernel in Theorem 2.4 of Lu--Tian\cite{LuTian2004} refers to the holomorphic extension of the
Szeg\H{o} kernel to the disk bundle, whereas the Bergman kernel considered here, as in \cite{EbenfeltXiaoXu2026}, is defined using holomorphic top-degree forms. Thus, our example does not contradict their theorem.

Let
\[
M=\operatorname{Bl}_p\mathbb{CP}^n,\qquad p=[1\colon 0\colon\cdots\colon 0].
\]
Realize $M$ as
\[
M=\bigl\{([Z_0\colon\cdots\colon Z_n],[u_1\colon\cdots\colon u_n])\in\mathbb{CP}^n\times\mathbb{CP}^{n-1}\colon Z_i u_j=Z_j u_i\quad(1\leq i,j\leq n)\bigr\},
\]
and let $L=\mathcal O(1,1)|_M$ be equipped with the Hermitian metric $h$ induced by the product Fubini--Study metric. Let $\pi_1\colon M\to\mathbb{CP}^n$ and $\pi_2\colon M\to\mathbb{CP}^{n-1}$ be the restrictions of the two projections. Denote by $H$ the pullback under $\pi_1$ of the hyperplane class and by $E=\pi_1^{-1}(p)$ the exceptional divisor. Note that
\[
\pi_1^*\mathcal O_{\mathbb{CP}^n}(1)=\mathcal O_M(H),\qquad\pi_2^*\mathcal O_{\mathbb{CP}^{n-1}}(1)=\mathcal O_M(H-E).
\]
Consequently,
\[
\begin{aligned}
L&=\pi_1^*\mathcal O_{\mathbb{CP}^n}(1)\otimes\pi_2^*\mathcal O_{\mathbb{CP}^{n-1}}(1)\\
&=\mathcal O_M(2H-E).
\end{aligned}
\]
By the blow-up formula for the canonical bundle,
\[
\begin{aligned}
K_M
&=\pi_1^*K_{\mathbb{CP}^n}\otimes\mathcal O_M((n-1)E)\\
&=\mathcal O_M(-(n+1)H+(n-1)E).
\end{aligned}
\]

On the dense open subset
\[
U=\pi_1^{-1}(\{Z_0\neq0\})\setminus E\cong\mathbb C^n\setminus\{0\},
\]
put
\[
z_j=\frac{Z_j}{Z_0},\qquad x=\sum_{j=1}^n|z_j|^2.
\]
Choose a holomorphic frame $e$ of $L$ over $U$ such that
\[
|e|_h^2=\frac1{x(1+x)}.
\]
Then the local K\"ahler potential is
\[
\phi(z)=\log x+\log(1+x).
\]
We use the convention $dV_g=\frac{\omega_g^n}{n!}$, denote the unit circle bundle of $(L^*,h^*)$ by $S$, and let $\pi\colon S\to M$ be the bundle projection. Throughout this appendix, the Szeg\H{o} kernel is defined with respect to $d\sigma=\pi^*dV_g\,\frac{d\theta}{2\pi}$.

A direct computation gives
\[
\begin{aligned}
g_{i\bar j}&=\phi'(x)\delta_{ij}+\phi''(x)\bar z_i z_j\\
&=\left(\frac1x+\frac1{1+x}\right)\delta_{ij}-\left(\frac1{x^2}+\frac1{(1+x)^2}\right)\bar z_i z_j,
\end{aligned}
\]
and hence
\[
\det g=(\phi')^{n-1}(\phi'+x\phi'')=\frac{(1+2x)^{n-1}}{x^{n-1}(1+x)^{n+1}}.
\]
The function $x$ extends smoothly to the exceptional divisor $E$, where it vanishes identically.

First, we consider the canonically twisted Bergman kernel. We have
\[
L^k\otimes K_M=\mathcal O_M\bigl((2k-n-1)H-(k-n+1)E\bigr).
\]
For $k\geq n$,
\[
H^0(M,L^k\otimes K_M)\cong H^0\bigl(\mathbb{CP}^n,\mathcal I_p^{k-n+1}(2k-n-1)\bigr).
\]
Thus the holomorphic sections correspond to homogeneous polynomials of degree $2k-n-1$ vanishing to order at least $k-n+1$ at $p$. Writing
\[
\alpha=(\alpha_1,\ldots,\alpha_n)\in\mathbb Z_{\geq0}^n,\qquad |\alpha|=\sum_{i=1}^n\alpha_i,\qquad \alpha!=\prod_{i=1}^n\alpha_i!,
\]
a basis is represented by
\[
Z_0^{2k-n-1-|\alpha|}Z_1^{\alpha_1}\cdots Z_n^{\alpha_n},\qquad k-n+1\leq|\alpha|\leq2k-n-1.
\]
On $U$, the corresponding sections are
\[
s_\alpha=z^\alpha e^{\otimes k}\otimes dz_1\wedge\cdots\wedge dz_n,\qquad k-n+1\leq|\alpha|\leq2k-n-1.
\]

The metric induced on $K_M$ satisfies
\[
|dz_1\wedge\cdots\wedge dz_n|_{K_M}^2=(\det g)^{-1},\qquad dV_g=2^n\det g\,dA(z),
\]
where $dA(z)=\prod_{i=1}^n dA(z_i)$ and
$dA(z_i)=\frac{\sqrt{-1}}2dz_i\wedge d\bar z_i$
denotes the Euclidean area form on $\mathbb C$.
Let $d=|\alpha|$. Since the complement of $U$ has measure zero, radial integration gives
\[
\begin{aligned}
\|s_\alpha\|_k^2&=2^n\int_{\mathbb C^n}|z^\alpha|^2e^{-k\phi(z)}\,dA(z)\\
&=\frac{(2\pi)^n\alpha!}{(d+n-1)!}\int_0^\infty x^{d+n-1-k}(1+x)^{-k}\,dx\\
&=\frac{(2\pi)^n\alpha!}{(d+n-1)!}\frac{\Gamma(d+n-k)\Gamma(2k-d-n)}{\Gamma(k)}.
\end{aligned}
\]
These sections are mutually orthogonal by invariance under coordinatewise rotations. Using
\[
\sum_{|\alpha|=d}\frac{|z^\alpha|^2}{\alpha!}
=\frac{x^d}{d!},
\]
we obtain
\[
B_k^{K_M}(z)=\frac{e^{-k\phi(z)}}{(2\pi)^n\det g}\sum_{d=k-n+1}^{2k-n-1}\frac{(d+n-1)!\,\Gamma(k)}{d!\,\Gamma(d+n-k)\Gamma(2k-d-n)}x^d.
\]

To evaluate this sum, write
\[
(t)_{\underline l}=t(t-1)\cdots(t-l+1),\qquad (t)_{\underline0}=1,\qquad y=\frac{x}{1+x}.
\]
Substituting $d=k-n+1+j$ gives
\[
B_k^{K_M}(z)=\frac{k-1}{(2\pi)^n}\frac{(1+x)^{n+1-k}}{(1+2x)^{n-1}}\sum_{j=0}^{k-2}\binom{k-2}{j}(k+j)_{\underline{n-1}}x^j.
\]
The identities
\[
(k+j)_{\underline{n-1}}=\sum_{l=0}^{n-1}\binom{n-1}{l}(k)_{\underline{n-1-l}}(j)_{\underline l}
\]
and
\[
\sum_{j=0}^{k-2}\binom{k-2}{j}(j)_{\underline l}y^j(1-y)^{k-2-j}=(k-2)_{\underline l}y^l
\]
yield the exact formula
\[
B_k^{K_M}(z)=\frac{k-1}{(2\pi)^n(1+y)^{n-1}}\sum_{l=0}^{n-1}\binom{n-1}{l}(k)_{\underline{n-1-l}}(k-2)_{\underline l}y^l.
\]
The right-hand side is a polynomial in $k$ of degree $n$, with leading coefficient $(2\pi)^{-n}$. Moreover, $y=\frac{\sum_{i=1}^n|Z_i|^2}
{\sum_{i=0}^n|Z_i|^2}$ pulls back to a smooth function on $M$. Thus the identity extends to all of $M$, and the canonically twisted Tian--Yau--Zelditch coefficients satisfy
\[
a_j^{K_M}\equiv0,\qquad j\geq n+1.
\]
By the characterization of Ebenfelt--Xiao--Xu \cite{EbenfeltXiaoXu2026}, the associated unit circle bundle is Bergman logarithmically flat.

Next, we consider the ordinary Bergman kernel. Since
\[
L^k=\mathcal O_M(2kH-kE),\qquad H^0(M,L^k)\cong H^0(\mathbb{CP}^n,\mathcal I_p^k(2k)),
\]
a basis is given on $U$ by
\[
t_\alpha=z^\alpha e^{\otimes k},\qquad k\leq|\alpha|\leq2k.
\]
Let $d=|\alpha|=k+j$, where $0\leq j\leq k$. The squared norm of $t_\alpha$ is
\[
\begin{aligned}
\|t_\alpha\|_k^2&=2^n\int_{\mathbb C^n}|z^\alpha|^2e^{-k\phi(z)}\det g\,dA(z)\\
&=\frac{(2\pi)^n\alpha!}{(k+j+n-1)!}\int_0^\infty\frac{x^j(1+2x)^{n-1}}{(1+x)^{k+n+1}}\,dx.
\end{aligned}
\]
Set
\[
I_{k,j}^{(n)}=\int_0^\infty\frac{x^j(1+2x)^{n-1}}{(1+x)^{k+n+1}}\,dx.
\]
Expanding $(1+2x)^{n-1}$ and evaluating the beta integrals gives
\[
I_{k,j}^{(n)}=\frac1{\Gamma(k+n+1)}\sum_{l=0}^{n-1}\binom{n-1}{l}2^l\Gamma(j+l+1)\Gamma(k+n-j-l).
\]
The sections are again mutually orthogonal. Consequently,
\[
B_k(z)=\frac1{(2\pi)^n(1+x)^k}\sum_{j=0}^k\frac{(k+j+n-1)!}{(k+j)!\,I_{k,j}^{(n)}}x^j.
\]
Along the exceptional divisor $E$, we have $x=0$, so only the term $j=0$ remains, and
\[
B_k|_E=\frac1{(2\pi)^n}\frac{(k+1)\cdots(k+n-1)}{I_{k,0}^{(n)}}.
\]
The substitution $t=x/(1+x)$ gives
\[
\begin{aligned}
I_{k,0}^{(n)}
&=\int_0^1(1+t)^{n-1}(1-t)^k\,dt\\
&=\sum_{l=0}^{n-1}\binom{n-1}{l}\frac{l!\,k!}{(k+l+1)!}\\
&=\frac{Q_{n-1}(k)}{\prod_{m=1}^n(k+m)},
\end{aligned}
\]
where
\[
Q_{n-1}(k)=\sum_{l=0}^{n-1}\binom{n-1}{l}l!\prod_{m=l+2}^n(k+m),
\]
hence
\[
B_k|_E=\frac1{(2\pi)^n}\frac{\displaystyle\prod_{m=1}^{n-1}(k+m)\prod_{m=1}^{n}(k+m)}{Q_{n-1}(k)}.
\]

The polynomial $Q_{n-1}$ is monic of degree $n-1$. For every integer $1\leq s\leq n$,
\[
Q_{n-1}(-s)=\sum_{l=s-1}^{n-1}\binom{n-1}{l}l!\prod_{m=l+2}^n(m-s)>0.
\]
Since the zeros of the numerator are contained in $\{-1,\ldots,-n\}$, no factor of $Q_{n-1}$ cancels. For $n\geq2$, the denominator is nonconstant, so $B_k|_E$ is not a polynomial in $k$.

It follows that its Laurent expansion at infinity contains a nonzero negative-power term. By uniqueness of the asymptotic expansion, comparison with
\[
B_k|_E\sim\left(\frac{k}{2\pi}\right)^n\sum_{j=0}^\infty\frac{a_j|_E}{k^j}
\]
gives
\[
a_j|_E\neq0 \qquad\text{for at least one }j\geq n+1.
\]
By the necessary condition of Lu--Tian \cite[Theorem~2.3]{LuTian2004} (see also \cite[Theorem~B]{PalmieriShen2026}), the unit circle bundle is not Szeg\H{o} logarithmically flat with respect to the natural measure.





%



\bibliographystyle{alpha}
\bibliography{ref.bib}

\end{document}